\documentclass[11pt,a4paper]{article}

\usepackage[english]{babel}
\usepackage[utf8]{inputenc}

\usepackage{amssymb,amsmath,amsthm,amsfonts}
\usepackage{txfonts}
\usepackage{graphicx}
\usepackage{verbatim}
\usepackage{color}
\usepackage[svgnames]{xcolor}
\usepackage{natbib}
\usepackage{enumitem} 
\usepackage{multicol}
\usepackage{booktabs,multirow}
\usepackage{tabularx}
\usepackage{subfloat}
\usepackage{subfig}
\usepackage{yhmath}
\usepackage{placeins}

\bibpunct{(}{)}{,}{a}{,}{,}

\usepackage[pdftex,citecolor=black,linkcolor=black,urlcolor=black,colorlinks=false]{hyperref} 

\newcommand{\ESS}{\text{ESS}}

\usepackage{fullpage}

\graphicspath{{RMSE/}{Coal/}{LRT/}{VraisProfileSimu/}{./}}

\usepackage{authblk}

\title{Resampling: an improvement of Importance Sampling in varying population size models}
\author[1,2,3,*]{C. Merle }
\author[2,3]{R. Leblois}
\author[3,4]{F. Rousset}
\author[1,3,5,$\dag$]{P. Pudlo }

\affil[1]{Institut Montpelli\'erain Alexander Grothendieck UMR CNRS 5149,  Universit\'e de Montpellier, Place Eug\`ene Bataillon, 34095 Montpellier, France }
\affil[2]{Centre de Biologie pour la Gestion des Populations, Inra, CBGP (UMR INRA / IRD / Cirad / Montpellier SupAgro), Campus International de Baillarguet, 34988 Montferrier sur Lez, France}
\affil[3]{Institut de Biologie Computationnelle, Universit\'e de Montpellier, 95 rue de la Galera, 34095 Montpellier}
\affil[4]{Institut des Sciences de l'Evolution, Universit\'e de Montpellier, CNRS, IRD, EPHE,  Place Eug\`ene Bataillon, 34392 Montpellier, France}
\affil[5]{Institut de Math\'ematiques de Marseille, I2M UMR 7373,  Aix-Marseille Universit\'e, CNRS, Centrale Marseille,
39, rue F. Joliot Curie, 13453 Marseille, France}
\affil[*]{\href{mailto:coralie.merle@umontpellier.fr}{\tt coralie.merle@umontpellier.fr}}
\affil[$\dag$]{\href{mailto:pierre.pudlo@univ-amu.fr}{\tt pierre.pudlo@univ-amu.fr}}

\date{}

\begin{document}
 \maketitle













\begin{abstract}
Sequential importance sampling algorithms have been defined to estimate likelihoods in models of ancestral population processes. However, these algorithms are based on features of the models with constant population size, and become inefficient when the population size varies in time, making likelihood-based inferences difficult in many demographic situations. In this work, we modify a previous sequential importance sampling algorithm to improve the efficiency of the likelihood estimation. Our procedure is still based on features of the model with constant size, but uses a resampling technique with a new resampling probability distribution depending on the pairwise composite likelihood.  We tested our algorithm, called sequential importance sampling with resampling (SISR) on simulated data sets under different demographic cases. In most cases, we divided the computational cost by two for the same accuracy of inference, in some cases even by one hundred. This study provides the first assessment of the impact of such resampling techniques on parameter inference using sequential importance sampling, and extends the range of situations where likelihood inferences can be easily performed.
\end{abstract}


\section{Introduction}

Under genetic neutrality, the distribution of the genetic polymorphism in a sample of individuals depends on
the evolution of the population size through unobserved stochastic processes. Typically, these stochastic
processes describe the evolution of the alleles at a given locus of the individuals from a population backward
to their Most Recent Common Ancestor (MRCA). When the population size is constant and finite, Wright-Fisher
models describe this evolution and the coalescent theory approximates these models when the population size is
large. In this context, the history (genealogy with mutations) of the observed sample is a latent
process. 
One of the major challenge to conduct a parametric inference analysis with these models is computing the
likelihood of the data at any point $\phi$ of the parametric space.  Indeed, the likelihood at $\phi$ is the
integral of the probabilities of each possible realization of the latent process. In population genetics, the
likelihood of an observed sample is an integral over the distribution of ancestral histories that may have led
to this sample. In this work, we consider a class of Monte Carlo methods based on Sequential Important
Sampling (SIS) which provides an estimate of the integral. In this scheme, the important sampling distribution
proposes paths of the process among those who contribute the most to the sum defining the likelihood.

For models of panmictic population with constant size, \cite{griffiths1994sampling} described an algorithm
wherein a proposal distribution suggests histories of the sample by stepwise reduction of the data set, either
by coalescence of two identical genes or by removal of a mutation on a single gene lineage. \citet[][Theorem
1]{stephens2000inference} characterized the optimal proposal distribution for a large class of time
homogeneous models, but not in the cases treated in this paper, that is varying population size. However, in
most cases (even in time-homogeneous models) the optimal distribution cannot be practically computed and has to
be approximated. \cite{de2004importance1} developed a method for constructing such approximations for any
model where the mutation process can be described as a Markov chain on gene types and \cite{de2004importance2}
extended this to subdivided population models. These methods have been further elaborated for stepwise
mutation models in a subdivided population by \cite{de2005stepwise}. 
The latent process, namely the history of the data, of models with past changes in population size exhibits
inhomogeneity in time.  Thus the previous theoretical arguments which derive efficient proposal distributions
in the literature are no longer applicable in this context. Nevertheless, as shown by
\citet{leblois2014maximum}, we can adapt an importance distribution from the importance distributions of
models with constant size populations. But their simulation tests demonstrated some limits of the algorithm,
which most importantly are large computation times for demographic scenarios with strong changes due to a Monte Carlo estimate of the likelihood with high variance.

Our aim in this work is to improve the accuracy of the likelihood estimation for the same
computational cost. One direction could be to derive a new importance sampling proposal distribution
like \cite{hobolth2008importance} did for the infinite site
model.  In this paper we chose another direction which consists of resampling among the paths
proposed by the importance distribution.  Indeed, the major difference between the constant
demographic scenario and size varying models is that the likelihood of the latter is an integral over a space
of much larger dimension, since we cannot remove the (random) times between events in the latent
process (ie. in the past history) from the integral which defines the likelihood. By contrast, in
constant demographic scenario, the likelihood is a sum of products of probabilities of each event
whatever their time because of the time homogeneity of the latent process.  Thus, even if the
importance distribution adapted by \citet{leblois2014maximum} were the most efficient distribution
among a certain class of distributions, it remains that the integral we estimate with importance
sampling in the general demographic case is an integral over a space of much larger dimension and
that the efficiency of importance sampling decreases with the dimension of the integral.

We took the opportunity of the paper to present rigorously the stochastic model with general demographic
scenario in Section~\ref{StochasticModelLikelihood}. In particular, Section~\ref{sub:Markovian} ends with a
writing of the likelihood as an integral in Eq.~\eqref{eq:likelihood} where one clearly sees that the integral
is of much larger dimension than in the constant demographic case because we cannot
remove the integrals over the random times, as explained in Section~\ref{sub:practical} after presenting the
sequential importance sampling (SIS) algorithm. The time  inhomogeneous latent process is part of the folklore in the
neutral population genetic literature, but has never been written down explicitly.

The major contribution of the paper is the addition and the calibration of a resampling procedure of
Section~\ref{Resampling} in the SIS algorithm, based on \cite{liu2001theoretical} and \cite{liu2008monte}. The
novelty is mostly in the choice of the resampling distribution that we propose in
Section~\ref{sub:PCLresampling}, which depends on both the current weight of the latent path and the pairwise
composite likelihood (PCL) of the current state of the latent process.  Section~\ref{Results} presents
numerical results on the likelihood estimates of simulated data sets. These results highlight the benefit due
to the proposed resampling distributions in the likelihood estimates. We then plug the likelihood estimates in
an inference method presented in Section~\ref{sub:inference}. The remainder of Section~\ref{sec:improvments}
highlights how our proposals improve the estimate of the parameter, and the likelihood surface around the
maximum likelihood estimate to compute confidence intervals (CIs).  We can thus confirm that the gain due to
resampling also benefits to the demographic parameter estimation. We end Section~\ref{sec:improvments} with a
discussion on cases where the data do not hold much information regarding the parameter of interest, leading
to flat likelihood surfaces.  Finally we show the relevance of our methods by presenting numerical results on a bat data set where strong
evidence for population contraction had been already provided by \citet{storz2002testing}. All computatiosn for this work were performed using
an updated version of the Migraine software \citep{rousset2007likelihood,rousset2012likelihood,leblois2014maximum}.

\section{The stochastic model and its likelihood}\label{StochasticModelLikelihood}

To illustrate our method we consider genetic data from individuals of a single population sampled at time $t=0$. Let $N(t)$ be the population size, expressed in number of genes, $t$ generations away from the sampling time $t=0$. We assume that $N(t)$ is a parametric function of $t$, see Section~\ref{UnderlyingDemographicModel} for examples. In this Section we focus only on data from a given locus.

\subsection{Stochastic model} \label{StochasticModel}

\cite{kingman1982coalescent}'s coalescent process is the usual model to describe ancestral relationships between gene copies of the sample under neutrality in a population of constant, but relatively large size.  We superimpose a mutation model on the coalescent process to describe genes modifications along lineages. Since the evolution is neutral, the coalescent is independent of the mutation process.  To describe the resulting process, we introduce a random vector $\mathbf H_t$, indexed by the set of possible types of genes (possible alleles) $E$: if $A\in E$, the component $\mathbf{H}_t(A)$ counts the number of genes of type $A$ at time $t$ (i.e., $t$ generations away, backward in time, from the sampling time) in the genealogy of the sample. The likelihood of the genetic data is given by the distribution of $\mathbf{H}_0$, which cannot be written as an explicit function of the parameter $\phi$.

\subsection{Markovian description of the evolution}
\label{sub:Markovian}
Actually, we only have at our disposal the following description of the process $\mathbf H_t$ forward in time. 

 Let $\mathbf{e}_A$ denote the vector indexed by $E$ whose components are all equal to $0$, except the $A$-component which is equal to one.  On one hand, the probability of a new lineage of type $A$ in the genealogy at time $t-\delta$, knowing that $\mathbf{h}$ is the value of the process at time $t$, is
\begin{align}
  \mathbb{P} \Big(\mathbf{H}_{t-\delta} = \mathbf{h} + \mathbf{e}_A\Big|\mathbf{H}_t = \mathbf{h} \Big)
  & = \frac{(|\mathbf{h}|+1)\,\mathbf{h}(A)}{2\, N(t)} \delta + o(\delta)
    \label{eq:ratecoal}
\end{align}
where $|\mathbf{h}| = \sum_{A\in E} \mathbf{h}(A)$ is the total number of lineages at time $t$ in the genealogy and $o(\delta)$ a quantity that is negligible in front of $\delta$ when $\delta\to 0$.
On the other hand the probability of a mutation of a gene of type $B$ at time $t$ to a gene of type $A$ at time $t-\delta$
, knowing that $\mathbf{h}$ is the value of the process at time $t$, is
\begin{align}
  \mathbb{P} \Big(\mathbf{H}_{t-\delta} = \mathbf{h} + \mathbf{e}_A-\mathbf{e}_B \Big|\mathbf{H}_t = \mathbf{h} \Big)
  & = \mu \, \mathbf{h}(B) \, p_{B,\, A}\, \delta  + o(\delta)
    \label{eq:ratemut}
\end{align}
where $p_{B,A}$ is the mutation probability from allele $B$ to allele $A$ forward in time, $\mu$ is the mutation rate per generation per lineage, and $o(\delta)$ is a quantity that is negligible in front of $\delta$ when $\delta\to 0$.  Additionally, when there is only one lineage in the genealogy, the distribution of the gene type (the allele) is supposed to be the stationary distribution $\psi(\cdot)$ of the transition matrix $p=\{p_{B,A}; B,A\in E\}$:
\[
\mathbb{P}\Big( \mathbf{H}_{t}=\mathbf{e}_A \Big| |\mathbf{H}_t|=1 \Big) = \psi(A).
\]
Hence, forward in time (i.e., when $t$ decreases), the coalescent based model $\mathbf{H}_t$ is a pure jump, continuous time Markov process taking values in the set of integers vectors indexed by $E$, namely $\mathbb{N}^E$. But the process is time inhomogeneous because the coalescence rate of Eq.~\eqref{eq:ratecoal} depends on the current date $t$ through the function $N(t)$. Note that Eq.~\eqref{eq:ratecoal}  and \eqref{eq:ratemut} can both be written
\[
\mathbb{P} \Big(\mathbf{H}_{t-\delta} = \mathbf{h}'\Big|\mathbf{H}_t = \mathbf{h} \Big)
   = \Lambda_t(\mathbf{h}'| \mathbf{h}) \delta + o(\delta)
\]
for any $\mathbf{h}\neq \mathbf{h}'$ in $\mathbb{N}^E$ where
\begin{equation}
\Lambda_t(\mathbf{h}'|\mathbf{h}) = \begin{cases}
  {(|\mathbf{h}|+1)\,\mathbf{h}(A)}\Big/\big({2\, N(t)}\big) & \text{if } \mathbf{h}'= \mathbf{h}+\mathbf{e}_A
  \\
  \mu \, \mathbf{h}(B) \, p_{B,\, A} & \text{if } \mathbf{h}'= \mathbf{h} + \mathbf{e}_A-\mathbf{e}_B
  \\
  0 & \text{otherwise}
\end{cases}
\label{eq:Lambda}
\end{equation}
defines the intensity matrix of the Markov process.

\bigskip

We can give a more explicit description of the process since it is a pure jump process, or, in other words, $\mathbf{H}_t$ is a piecewise constant function of $t$.  We denote $\mathbf{X}_0=\mathbf{H}_0$ the first value of the process at time $T_0=0$. The process $\mathbf{H}_t$ remains constant until (random) time $T_1$, where it takes another value $\mathbf{X}_1 = \mathbf{H}_{T_1}$. After $T_1$, the process $\mathbf{H}_t$ stays equal to $\mathbf{X}_1$ until time $T_2$ where it jumps to a another value $\mathbf{X}_2$,  and so on. 

Fig.~\ref{fig:coal} in the Supplementary Material represents a possible path of the process $\mathbf{H}$. We refer to such paths as (possible) histories of the sample (composed of three genes of type a and one gene of type b in Fig.~1 in the supplementary material). Note that many genealogies correspond to a possible history since, at time $T_1$ we have chosen a genealogy with a coalescence between the two left-handmost genes, but we could have chosen another genealogy joining any pair of genes of type $a$ at time $T_1$. Likewise at time $T_3$ we could have chosen any of the three possible pairs. Actually each genealogy leading to the same path of the process $\mathbf{H}$ has the same probability because of the exchangeability of the genes carrying the same allele.

\bigskip

We set $\Delta_i = T_{i+1}-T_{i}$ which is usually named the holding time at value $\mathbf{X}_i$ for any integer number $i$. The distribution of the process is given by the density
\begin{multline}
  \mathbb{P}\Big(\mathbf{X}_{i}=\mathbf{h}_{i} \text{ for all }i=0,\ldots,n\text{ and }
  \Delta_i\in (\delta_i; \delta_i+d\delta_i) \text{ for all }i=0,\ldots,n-1 \Big|
\mathbf{X}_n=\mathbf{h}_n\Big)=\\
\prod_{i=1}^{n} P_{t_i}(\mathbf{h}_{i-1}| \mathbf{h}_{i}) \times
\prod_{i=0}^{n-1} \lambda_{t_i+\delta_i}(\mathbf{h}_i)\exp\left(-\int_{0}^{\delta_i}\lambda_{t_i+u}(\mathbf{h}_i)du\right)\,d\delta_i
\label{eq:density}
\end{multline}
for any $\delta_0,\ldots,\delta_{n-1}>0$, $\mathbf{h}_0, \ldots, \mathbf{h}_n$ in $\mathbb{N}^E$,
where $t_0=0$, $ t_i  = \delta_0 + \cdots + \delta_{i-1}$,
\begin{equation}
\lambda_t(\mathbf{h})  = \sum_{\mathbf{h}'\neq \mathbf{h}} \Lambda_t(\mathbf{h}'| \mathbf{h})
\label{eq:lambda}
\end{equation}
\[ \text{ and } \quad
  P_t(\mathbf{h}'| \mathbf{h}) = \Lambda_t(\mathbf{h}'| \mathbf{h}) \Big/ \lambda_t(\mathbf{h}).
\]
The term $\lambda_t(\mathbf{h})$ is interpreted as the (infinitesimal) jump rate of the process $\mathbf{H}$ at time $t$, knowing that the process takes value $\mathbf{h}$ at time $t$, namely
\[
\mathbb{P}\Big(\mathbf{H}_{t-\delta}\neq \mathbf{h}\Big|
\mathbf{H}_t=\mathbf{h} \Big) = \lambda_t(\mathbf{h})\,\delta + o(\delta)
\]
when $\delta>0$ and $o(\delta)$ is negligible in front $\delta$ when $\delta\to0$.
With Eq.~\eqref{eq:Lambda}, it can be computed as
\[
\lambda_t(\mathbf{h}) = \frac{|\mathbf{h}|(|\mathbf{h}|+1)}{2N(t)} + \mu |\mathbf{h}| =
\frac{|\mathbf{h}|}{2N(t)} \left( (|\mathbf{h}|+1) + \theta(t) \right),
\]
where $\theta(t) = 2\mu N(t)$ is the usual composite mutation rate parameter of population genetics.
On the other hand, $P_t(\mathbf{h}'|\mathbf{h})$ is interpreted as a transition matrix (forward in time), and gives the probability that the new value (forward in time) of $\mathbf{H}_{t'}$ is $\mathbf{h}'$ for $t'<t$ knowing that the process jumps at time $t$ and that $\mathbf{H}_t=\mathbf{h}$.

\bigskip

Set $\tau = \inf\big\{t>0:\ |\mathbf{H}_t|=1\big\}$, which is the age of the most recent common ancestor (MRCA) of the sample. Since $\mathbf{H}$ is a pure jump process, $\tau = T_\sigma$
where $\sigma = \inf\big\{n>0:\ |\mathbf{X}_n|=1\big\}$.
And note that the distribution of $\mathbf{H}_\tau=\mathbf{X}_\sigma$ is known explicitly as
\[
\mathbb{P}(\mathbf{H}_\tau = \mathbf{e}_A) = \psi(A).
\]
To recover the likelihood of the observed genetic data, that is to say the distribution of $\mathbf{H}_0$ from this knowledge and the (forward in time) transition mechanism of Eq.~\eqref{eq:density}, we have to integrate over all possible histories from the MRCA, that is to say all possible values $n$ of $\sigma$, $\mathbf{h}_1,\ldots, \mathbf{h}_n$ in $\mathbb{N}^E$ and $\delta_0,\ldots, \delta_{n-1}>0$.  Hence
\begin{multline}
\mathbb{P}\big(\mathbf{H}_0=\mathbf{h}_0\big)
  = \sum_{n=1}^\infty \int\cdots\int \sum_{\mathbf{h}_1, \cdots, \mathbf{h}_{n}}\sum_{A} \psi(A)\,
\mathbf 1\left\{\mathbf{h}_n = \mathbf{e}_{A}, |\mathbf{h}_{n-1}|>1\right\} \times
\\
\prod_{i=1}^{n} P_{\delta_0+\cdots+\delta_{i-1}}(\mathbf{h}_{i-1}| \mathbf{h}_{i}) \times
\prod_{i=0}^{n-1} \lambda_{\delta_0+\cdots+\delta_{i}}(\mathbf{h}_i)
\exp\left(-\int_{0}^{\delta_i}\lambda_{\delta_0+\cdots+\delta_{i-1}+u}(\mathbf{h}_i)du\right)\,d\delta_i.
\label{eq:likelihood}
\end{multline}
The challenge to conduct a likelihood based inference in population size varying models is in computing
the multidimensional integral of Eq.~\eqref{eq:likelihood}, which cannot be computed formally.
Note that, to alleviate notations, we have dropped the dependency in the parameter of interest $\phi$ of the
quantities arising in Eq.~\eqref{eq:likelihood} but both $P_t(\mathbf{h}'| \mathbf{h})$ and $\lambda_t(\mathbf{h})$ are functions of $\phi$, and even sometimes the stationary distribution $\psi(A)$.

\subsection{Evaluating the likelihood with importance sampling}

The trick to compute the likelihood defined in Eq.~\eqref{eq:likelihood} is to rely on another
transition matrix $Q_t(\mathbf{h}| \mathbf{h}')$ backward in time and set
\begin{equation}
W_n(\mathbf{h}_{0:n}; t_{0:n})=
{\displaystyle \prod_{i=1}^{n} P_{t_i}(\mathbf{h}_{i-1}| \mathbf{h}_{i})} \bigg/ {
\displaystyle\prod_{i=1}^{n} Q_{t_i}(\mathbf{h}_{i}| \mathbf{h}_{i-1})}.
\label{eq:Wn}
\end{equation}
Then Eq.~\eqref{eq:likelihood} leads to
\begin{multline}
\mathbb{P}\big(\mathbf{H}_0=\mathbf{h}_0\big)
  = \sum_{n=1}^\infty \int\cdots\int \sum_{\mathbf{h}_1, \cdots, \mathbf{h}_{n}}\sum_{A} \psi(A)
\mathbf 1\left\{\mathbf{h}_n = \mathbf{e}_{A}, |\mathbf{h}_{n-1}|>1\right\} \times
\\
W_n(\mathbf{h}_{0:n};  t_{0:n}) \times
\prod_{i=1}^{n} Q_{t_i}(\mathbf{h}_{i}| \mathbf{h}_{i-1}) \times
\prod_{i=0}^{n-1} \lambda_{t_i+\delta_i}(\mathbf{h}_i)\exp\left(-\int_{0}^{\delta_i}\lambda_{t_i+u}(\mathbf{h}_i)du\right)\,d\delta_i,
\label{eq:IS}
\end{multline}
where $t_i=\delta_0+\cdots+\delta_{i-1}$ are implicit functions of the $\delta_i$'s.
The right hand side of Eq.~\eqref{eq:IS} can be interpreted as the expected value of
\[
W_\sigma =
\sum_n\sum_{A} \psi(A)
\mathbf 1\left\{\mathbf{h}_n = \mathbf{e}_{A}, |\mathbf{h}_{n-1}|>1\right\} \times
W_n(\mathbf{h}_{0:n}; t_{0:n})
\]
when $\mathbf{h}_0, \mathbf{h}_1,\ldots$ and $t_0, t_1,\ldots$ are the realizations of the embedded chain and the jump times of an inhomogeneous Markov process $\widetilde{\mathbf{H}}$ whose intensity matrix (backward in time) is given by
\[
\widetilde\Lambda_t(\mathbf{h}| \mathbf{h}') = \lambda_t(\mathbf{h}')Q_t(\mathbf{h}| \mathbf{h}')
\]
and which starts from $\widetilde{H}_0=\mathbf{h}_0$, the observed data. In other words,
\[
\mathbb{P}(\mathbf{H}_0=\mathbf{h}_0) = \widetilde{\mathbb{E}}(W_\sigma).
\]
Interpreting the likelihood of the data as an expected value over another distribution is the first step toward an importance sampling estimate of the likelihood. Indeed, the Monte Carlo estimation of the likelihood is the empirical average of  $W_\sigma$ computed on simulated replicates $\widetilde{\mathbf{H}}^{(j)}$ (for $j=1,\ldots, n_H$) of the process $\widetilde{\mathbf{H}}$:
\begin{equation}
\widehat{\mathbb{P}}(\mathbf{H}_0=\mathbf{h}_0) = \frac{1}{n_H}\sum_{j=1}^{n_H} W_\sigma^{(j)}.
\label{eq:ISestimate}
\end{equation}

\subsection{Practical aspects and efficiency}
\label{sub:practical}
To sum up the above, the idea of importance sampling is to interpret the likelihood as an expected
value over a sampling distribution of histories (the distribution of the process $\widetilde{\mathbf{H}}$) which may differ from the distribution of the genuine model described in Section~\ref{sub:Markovian}. We can choose this importance distribution freely, as long as any path from the data $\mathbf{h_0}$ to the MRCA with positive density under the distribution of the latent process $\mathbf{H}$ has also a positive density under the distribution of the importance process $\widetilde{\mathbf{H}}$.  But this choice has a major impact on the efficiency of the approximation in Eq.~\eqref{eq:ISestimate} since the variance of $W_\sigma$ depends on the distribution of $\widetilde{\mathbf{H}}$.

\bigskip

When the population size $N(t)$ is fixed to $N_0$ for all $t$, we have at our disposal an efficient sampling distribution from the literature \citep{stephens2000inference, de2004importance1, de2005stepwise} which represents a major improvement over the first proposal of \cite{griffiths1994ancestral}. In this simple demographic scenario, the importance sampling estimate is the most efficient under a parent independent mutation model (which means that $p_{B A}$ does not depend on $B$) and Eq.~\eqref{eq:IS} provides an exact evaluation of the likelihood for $n_H=1$ replicate. In other words, the variance of $W_\sigma$ under the distribution of the efficient importance sampling distribution is zero. For other mutation models and a constant demographic scenario, the variance of $W_\sigma$ is no longer zero, but the number $n_H$ of replicates required to get a sharp estimation of the likelihood with Eq.~\eqref{eq:ISestimate} is much smaller when relying on the efficient importance distribution rather than on the proposal of \cite{griffiths1994ancestral}.
In the scenario where the population size varies over time,
we resort to the proposal of \cite{leblois2014maximum} to define the transition matrix of $\widetilde{\mathbf{H}}$. That is to say that for any value of $t>0$, $Q_t(\mathbf{h}, \mathbf{h}')$ is the transition matrix of the efficient importance distribution as if the population size were constant over time and fixed to the current value of $N(t)$.
But the efficiency of this importance distribution depends on the variation of $N(t)$, as explained in \cite{leblois2014maximum}.

\bigskip

Nevertheless, the choice of importance distribution leads to a process $\widetilde{\mathbf{H}}$ which is defined explicitly as an inhomogeneous Markov process, backward in time, starting from the observed data $\mathbf{h}_0$. The computation of a simulated path of $\widetilde{\mathbf{H}}$, as well as of $W_\sigma$ is performed with the sequential method of Algorithm~1.
\begin{table*}[hbt]
\hrule
\vspace{1mm}
\textbf{Algorithm~1}: sequential importance sampling
\vspace{1mm}
\hrule
\begin{enumerate}
	\item Initialization: set $\mathbf{h}=\mathbf{h}_0$, $w=1$ and $t=0$
	\item  \textbf{Repeat}
          \vspace*{-3mm}
	\begin{enumerate}[label=2.(\alph*)]
 	 	\item  Draw the holding time $\delta$ according to its density
                  $\displaystyle
                  \lambda_{t+\delta}(\mathbf{h})\exp\left(-\int_0^{\delta}\lambda_{t+\delta}(\mathbf{h})\right)
                  $\\
                  and update $t = t + \delta$
                \item  Draw $\mathbf{h}'$ according to $ q_{t}(\mathbf{h}'|\mathbf{h})$
                and  update $w = w \times p_{t}(\mathbf{h}|\mathbf{h}') / q_{t}(\mathbf{h}'|\mathbf{h})$
                \item  Update $\mathbf{h}= \mathbf{h}'$
	\end{enumerate}
        \hspace*{-.7cm} \textbf{until} $|\mathbf{h}|=1$
       \item  Set $W_\sigma = w \times \psi(\mathbf{h})$ and return $W_\sigma$
\end{enumerate}
\hrule
\end{table*}

Note that, since $W_\sigma$ is a product along the path of the process
$\widetilde{\mathbf{H}}$, its value is computed sequentially at step
2.(b) of Algorithm~1. The update depends only on the current time $t$,
the value of $\widetilde{\mathbf{H}}$ just before time $t$, namely
$\mathbf{h}$, and the new value $\mathbf{h}'$. Hence, Algorithm~1 does
not have to keep track of the whole path of $\widetilde{\mathbf{H}}$.
Moreover, drawing the holding time $\delta$ can be done either with a
rejection algorithm, see Appendix~\ref{AppendixGillespie} in the Supplementary Materials, or by inverting the cumulative
distribution function when possible, e.g., see \cite{griffiths1994sampling}.
Finally, the above algorithm is run $n_H$ times to approximate the
likelihood by the average in Eq.~\eqref{eq:ISestimate}.

\bigskip

Sadly, the decrease of efficiency of the importance sampling scheme from the constant
population size scenario to the general case (where $N(t)$ varies over time)
can be drastic, see \cite{leblois2014maximum}. To understand the major difference, we
recall that, when the population size $N(t)$ is constant, neither
$\Lambda_t(\mathbf{h}'| \mathbf{h})$ nor $P_t(\mathbf{h}'|\mathbf{h})$ depend on $t$. In
this simpler case we choose an importance distribution characterized by a
$Q_t(\mathbf{h}|\mathbf{h}')$ which does not depend on $t$, and $W_\sigma$ does not depend
on the random jump times $T_1,T_2,\ldots$ of the process $\widetilde{\mathbf{H}}$. But in
the varying size scenario, these random times do contribute to the variance of $W_\sigma$.

Fig.~\ref{0.4_0_0.25_400_NormaliseWeightTrajectory} follows $w$, computed at line
2.(b) of Algorithm~1, over successive coalescence for one hundred replicates of the SIS algorithm. It shows first that the variance of the final SIS weights is quite large, second that replicates that lead to the highest final $W_\sigma$ also tend to have high $w$ throughout the sequence of coalescence events.
Moreover, Fig.~\ref{0.4_0_0.25_400_NormaliseWeightTrajectory} shows that, as the number of coalescence events
undergone increases, the range of values of $w$ increases exponentially. This exponential increase
is an evidence that importance sampling becomes inefficient in spaces of high dimension, and that
each random jump times $T_1, T_2,\ldots$ has a multiplicative contribution to the overall range of
$W_\sigma$.

The final approximation of the likelihood is the average of these one hundred replicates
of $W_\sigma$, see Eq.~\eqref{eq:ISestimate}. The difference between the values of the
final SIS weights $W_\sigma$ is so huge that most of them are negligible compared to the
top ten highest ones and do not contribute to the average in Eq.~\eqref{eq:ISestimate}. If we look at Fig.~\ref{0.4_0_0.25_400_NormaliseWeightTrajectory},
considerations of their trajectories further back in time show that these weights do not
rise back much. Hence, we can reduce the computing effort and increase the accuracy of the
method by ignoring the histories which give low values of $w$ after some times through
the sequence of coalescence events and by replacing them by histories with higher values of SIS weight
$w$. The above idea can be seen as a correction of the proposal distribution, and the
estimate  will remain unbiased if the SIS weights are accordingly corrected. In order to
implement it, we rely on the resampling procedure explained in the next Section.

\section{Resampling}\label{Resampling}

We derive a resampling procedure from that of \citet[][Section~4.1.2]{liu2008monte} that
consists of pruning away the partial histories associated with very small weights
and reusing those associated with high weights instead of restarting from scratch.

\subsection{Sequential importance sampling with resampling: the algorithm}
\label{sub:SISR}

\begin{table*}[h!]
\hrule
\vspace{1mm}
\textbf{Algorithm~2}: Sequential importance sampling with resampling (SISR)
\vspace{1mm}
\hrule
\vspace*{-3mm}
\begin{enumerate}[label=\bf\Roman*.]
\item Initialization:\\
  \textbf{For} $j$ in $1,\ldots, n_H$:\\
  \hspace*{.8cm} Set $w^{(j)}=1$, $t^{(j)}=0$ and $\mathbf{h}^{(j)}=\mathbf{h}_0$\\
  \textbf{EndFor}\\
  Set $\text{ESS}_-=n_H$

\item  \textbf{Repeat}
  \begin{enumerate}[label={\bf II.} \arabic*]
  \item \textbf{For} $j=1, \ldots,  n_H$
    \begin{itemize}[label=\ \ \ ]
    \item Set $n_{control} = 0$
    \item
      \textbf{While}   $|\mathbf{h}^{(j)}| > 1$ \textbf{and} $n_\text{control} < k$:
      \begin{enumerate}[label=(\alph*)]
      \item  Draw the  holding time $\delta$ according to its distribution \\
        and update $t^{(j)} = t^{(j)} + \delta$
      \item  Draw $\mathbf{h}'$ according to $q_{t^{(j)}}(\mathbf{h}'|\mathbf{h}^{(j)}) $\\
        and update $w^{(j)} = w^{(j)} \times p_{t^{(j)}}(\mathbf{h}^{(j)}|\mathbf{h}') / q_{t^{(j)}}(\mathbf{h}'|\mathbf{h}^{(j)})$
      \item  Update $n_\text{control}$ according to its definition and update $\mathbf{h}^{(j)} = \mathbf{h}'$
      \end{enumerate}
      \textbf{EndWhile}
    \end{itemize}
    \textbf{EndFor}
    \vspace*{-3mm}
  \item Compute $\displaystyle \ESS_+ = \left(\sum_{j=1}^{n_H} w^{(j)}\right)^2\bigg/\sum_{j=1}^{n_H} {\left(w^{(j)}\right)^2}$
  \item \textbf{If} $\ESS_+ <  \ESS_-/10$, \textbf{then}
    \\
    \hspace*{.8cm} Resample according to Algorithm~3\\
    \hspace*{.8cm} Update $\ESS_-$ to the $\ESS$ of the resampled collection
    \\
    \textbf{EndIf}
  \end{enumerate}
  \hspace*{-.8cm}\textbf{Until} all $|\mathbf{h}^{(j)}|$ are equal to $1$
\item    \textbf{For} $j$ in $1,\ldots, n_H$:\\
  \hspace*{.8cm} Update $ w^{(j)} = w^{(j)} \times \psi(\mathbf{h}^{(j)})$\\
  \textbf{EndFor}
\item  Return the average $\displaystyle n_H^{-1}\sum_{j=1}^{n_H} w^{(j)}$
\end{enumerate}
\vspace*{-3mm}
\hrule
\vspace{1mm}
\end{table*}

The sequential importance sampling with resampling (SISR) is as follows. We initiate $n_H$ independent runs of Algorithm~1, hence $n_H$ draws of $W_\sigma$ from the importance distribution.
But we stop the repeat-until loop at step 2 of Algorithm~1 before hitting the MRCA. Indeed each run of Algorithm~1 is stopped at checkpoints (stopping times of the Markov process $\widetilde H_t$); at these checkpoints we evaluate the quality of the collection of $n_H$ partial histories and if necessary, we resample these histories. To that aim, we propose a new resampling distribution detailed in Section~\ref{sub:resampling}.

The checkpoints at which we test whether to resample correspond either to a given number $k$ of events (coalescences and mutations) undergone after the previous checkpoint or to a given number $k$ of coalescence events undergone.
Once a checkpoint is reached, we evaluate the quality of the collection of partial histories to assess whether it is necessary to resample. The quality test is based on the Effective Sample Size (ESS) relative to the collection of partial histories, namely
\begin{equation}
 \text{ESS} = { \left(\sum_{j=1}^{n_H} {w^{(j)}}\right)^2} \bigg/ {\sum_{j=1}^{n_H} {(w^{(j)})^2}}
\label{9}
\end{equation}
where $w^{(j)}$ is the current value of $w$ of the $j$-th partial history at the checkpoint. The largest value of ESS is $n_H$ and occurs when $w^{(1)}=\cdots=w^{(n_H)}$; the ESS decreases when the range of values of $w^{(j)}$ expands. When a single weight, $w^{(1)}$ say, is much larger than the other ones, the ESS is approximately equal to $1$ since both numerator and denominator of Eq.~\eqref{9} are approximately equal to $\big(w^{(1)}\big)^2$. Actually, the ESS assesses how the random holding times and events of the partial histories between two successive checkpoints contribute to the variance of the estimate in Eq.~\eqref{eq:ISestimate}.
Thus we resample the collection of partial histories whenever the ESS falls below a threshold value, for instance $\text{ESS}_-/10$, where $\text{ESS}_-$ is the value of the ESS after the last resampling.

The checkpoints are frequent when $k=1$, but computing the ESS, as well as resampling the whole collection is time consuming, so that higher values of $k$ might be more pragmatic.

The SISR is presented in Algorithm~2, where $w^{(j)}$, $t^{(j)}$ and $\mathbf{h}^{(j)}$, for $j=1,\ldots, n_H$ are arrays which keep tracks of the values of current weight $w$, time $t$ and state $\mathbf{h}$ of each run of Algorithm~1. Moreover, we define $n_\text{control}$ which stores the number of undergone events (either the total number of undergone events, or the number of undergone coalescences) inbetween checkpoints.

\bigskip

\subsection{The resampling procedure}
\label{sub:resampling}

Assume that SISR has reached a checkpoint and that the $\ESS$ is low enough to resample. At this time, we create a new collection of $n_H$ simulated histories by drawing at random in the previous collection of histories according to a multinomial distribution $\mathcal{M}(\mathbf{v},n_H)$ where $\mathbf{v}=(v^{(1)}, \ldots, v^{(n_H)})$ is a resampling probability distribution on the collection of $n_H$ partial histories, see Section~\ref{sub:PCLresampling} for examples of such distributions.

Resampling is equivalent to applying a second importance sampling algorithm within the SIS. Indeed,
\begin{equation}\label{eq:IS2}
\sum_{j=1}^{n_H} w^{(j)} = \sum_{j=1}^{n_H} \frac{w^{(j)}}{v_j^{(j)}} v^{(j)}
\end{equation}
and the right hand side can be interpreted as the expected value of $w^{(j)}/v^{(j)}$ over the distribution of a random $j$ drawn from $\mathbf{v}$.
Thus, in order not to bias the procedure, the new weight associated to the $j$-th history is $w^{(j)}/v^{(j)}$ whenever this history appears in the resampled collection of histories. Algorithm~3 summarizes the procedure, with the notations of Algorithm~2.

\begin{table*}[h!]
\hrule
\vspace{1mm}
\textbf{Algorithm~3} Resampling procedure
\vspace{1mm}
\hrule
\vspace*{-3mm}
\begin{enumerate}
\item \textbf{For} \( j = 1,...,n_\text{H},\)
\begin{enumerate}[label=(\alph*)]
\item Draw $J'$ from distribution $\mathbf{v}$
\item Set $\widetilde{\mathbf{h}}^{(j)}=\mathbf{h}^{(J')}$ and $\widetilde{t}^{(j)}=t^{(J')}$
\item Set $\widetilde{w}^{(j)}=w^{(J')} / v^{(J')}$
\end{enumerate}
\hspace*{-.7cm} \textbf{EndFor}
\item Replace the old collection \( \left\{( \mathbf{h}^{(j)}, t^{(j)}, w^{(j)})\right\}^{n_H}_{j=1}\)
with the new collection\\ \( \left\{(\widetilde{\mathbf{h}}^{(j)}, \widetilde{t}^{(j)}, \widetilde{w}^{(j)})\right\}^{n_H}_{j=1}\).
\end{enumerate}
\vspace*{-3mm}
\hrule
\end{table*}

\subsection{The resampling distribution}
\label{sub:PCLresampling}

Resampling introduces a new possible cause of variance, since we replace the sum of Eq.~\eqref{eq:IS2} with a Monte Carlo estimate.
The resampling probability distribution $\mathbf{v}=(v_1, \ldots, v_{n_H})$ on the collection of histories could be any distribution.
To achieve efficiency we should pick a distribution $\mathbf{v}$ that reflects the future trend of the partial histories, more precisely that sets relatively high probabilities on the partial histories that will correspond to the highest $W_\sigma$.
If the resampling distribution brings helpful information about histories, then the mean square error (MSE) of Eq.~\eqref{eq:ISestimate} should decrease. On the contrary if the resampling probability distribution brings useless information or no information, then the resampling just adds noise and the MSE increases.
For example a uniform $\mathbf{v}$ distribution that clearly does not bring any information about the histories introduces only an additional variance in the estimation of the likelihood.

\bigskip
The resampling algorithm could face two difficulties: first it does not always choose the replicates with the highest $W_\sigma$; second, regardless of the optimal replicate, it is not perfectly predicted from an intermediary $w$.
However, as explained at the end of Section~\ref{sub:practical}, the current value of $w$ helps predicting the contribution of the final value $W_\sigma$ to the average in Eq.~\eqref{eq:ISestimate}.
One might thus be tempted to resample the $j$-th partial history with a probability proportional to $w^{(j)}$, i.e. $v_j \propto w^{(j)}$. But this resampling distribution $\mathbf{v}$ might not well choose the optimal $W_\sigma$ because of the spread of the $w^{(j)}$'s. Thus
\citet{liu2001theoretical} proposed to rely on $v_j \propto \left[w^{(j)}\right]^\alpha$ for some $\alpha\in [0; 1]$, with an arbitrary preference on $\alpha = 1/2$.
Another possible predictor of the contribution of the $j$-th partial history to the average in Eq.~\eqref{eq:ISestimate} is its current state $\mathbf{h}^{(j)}$. Indeed, if the probability of $\mathbf{h}^{(j)}$ under the distribution of the process defined in Eq.~\eqref{eq:density} is large, then the final weight $W_\sigma$ tends to be high, while if its probability is small, $W_\sigma$  might be negligible.
Of course the probability of $\mathbf{h}^{(j)}$ is intractable, but we might replace it by an easily computed pseudo-likelihood. We propose here to rely on the pairwise composite likelihood $L_2(\mathbf{h}^{(j)})$ defined in the Appendix~\ref{AppendixPCL} in the Supplementary Materials since we can compute it very easily. But such pseudo-likelihoods are much more contrasted than the true one and thus should be tempered by some exponent $\beta \ll 1$.
Hence we advocate the following resampling distribution
\begin{equation}
 v^{(j)} \propto \left(w^{(j)}\right)^\alpha \left(L_2(\mathbf{h}^{(j)})\right)^\beta, \quad \alpha, \beta \in \left[0,1\right], \beta \ll 1.
\label{NewWeigths}
\end{equation}
The tuning parameters $\alpha$ and $\beta$ are used to balance the effect of the information provided by the SIS weight and by the composite likelihood. Section~\ref{Results} provides numerical examples showing the efficiency of the above resampling distribution for a large range of values of the tuning parameters $\alpha$ and $\beta$.

\section{Improvements on the likelihood estimate: numerical results} \label{Results}

\subsection{The simulated demographic model}  \label{UnderlyingDemographicModel}
We use the model described in Sections~\ref{StochasticModel} and  \ref{sub:Markovian} to analyze microsatellite markers under a Stepwise Mutational Model (SMM). The set $E$ of allele types corresponding to microsatellite markers is the set $\mathbb{N}$ of non negative integers. When a mutation occurs under a SMM, the size of the allele is either increased by $1$ or decreased by $1$ with the same probability \citep{kimura1978stepwise}.

As in \cite{leblois2014maximum} we consider a single isolated population whose size has undergone past changes. We denote \( {N}(t)\) the population size expressed as the number of genes, \(t\) generations away from the sampling time \(t=0\). The population size at sampling time is \( {N}(0)=N\).  Then, going backward in time, the population size changes according to a deterministic function until reaching an ancestral population size \(N_\text{anc}\) at time \(t = T\). Then, \( {N}(t) = N_\text{anc}\) for all \(t > T\).
To illustrate our method we consider an exponentially contracting population size, represented in Fig.~\ref{ECP} in the supplementary material, which can be written
\begin{equation*}
 {N}(t) =  \begin{cases}
 N \times \left(N_\text{anc}/N\right)^{t/T} & \text{if } t \in [ 0; T ], \\
N_\text{anc} & \text{if } t \geq T.
\end{cases}
\label{N(t)expo}
\end{equation*}

To ensure identifiability the demographic parameters are scaled as \(\theta=2 \mu N\), \(D=T/2N\) and \(\theta_\text{anc}= 2 \mu N_\text{anc}\), where \(\mu\) is the mutation rate per locus per generation. The parameter space of the model is thus the set of vectors $\phi=(\theta, D, \theta_\text{anc})$. Additionally we also set  \(\theta(t) =2 \mu  {N}(t) \), and we are sometimes interested in the extra composite parameter \(N_\text{ratio}=\theta / \theta_\text{anc}= N / N_\text{anc}\), which is useful to characterize the strength of the contraction.

\bigskip

\subsection{Reduction of the MSE between the true value of the likelihood and its estimate}
\label{LikEst}

Numerical results of Section~\ref{Results} are presented on data sets simulated under the demographic scenario where $\phi =(\theta, D,\theta_\text{anc})=(0.4, 0.25, 400)$ which models a recent and strong past contraction in population size. This scenario is the most representative of our conclusions although we studied more moderate and/or older contractions.
To assess the efficiency of estimates of the likelihood at a given point of the parameter space, we compared these estimates with a reference value $L$ that is a sharp approximation of the true likelihood.  We have computed the reference value $L$ by estimating the likelihood with SIS with a huge collection of $n_H=20,000,000$ independent histories.

To evaluate the variability of the likelihood estimate returned by a given algorithm, with given values of its tuning parameters, we first plot the empirical distribution of $100$ runs of the same algorithm with a boxplot (Fig.~\ref{Boxplotnevtncoal}--\ref{Boxplots}). We also measured the mean square error of $100$ estimates $\widehat{L_i}$ around the reference value $L$, and averaged this MSE over $100$ simulated data sets to remove the dependency on the simulated data set (Fig.~\ref{0.4_0_0.25_400_RMSE}).
Note that this last method computes an empirical mean square error that approximates the MSE of the likelihood estimate.

To compare the efficiency of both SIS and SISR algorithms with the same computational effort of simulating a collection of $n_H=100$ histories, we considered the ratio of the empirical mean square error of SISR over that of SIS. The MSE ratio brings out the reduction of the mean square error in likelihood estimate due to the resampling technique. 

\bigskip

For Fig.~\ref{Boxplotnevtncoal}--\ref{Boxplots}
 discussed below, we note that all the SISR calibrations lead to a smaller variation in likelihood estimation than SIS and that the average is much closer to the reference value meaning that resampling improve the accuracy of the likelihood estimate. This observation means that even if the resampling procedure is not much calibrated, the likelihood inference is already better than without resampling.

\paragraph{Type and frequency of checkpoints}
As said in Section~\ref{sub:SISR}, the checkpoints might be defined either in term of a fixed number of events (coalescences and mutations) undergone on each partial history or in term of a fixed number of coalescence events. Note that the total number of genes in the current state $\mathbf{h}$ is determined by the size of the observed data $\mathbf{h}_0$ and the number of coalescence events undergone from time $t=0$ to state $\mathbf{h}$. Hence the second definition of a checkpoint proposes to compare partial histories that lead to current states $\mathbf{h}$ with the same number of lineages. Thus \citet{liu2008monte} advocated the second type of checkpoints.
Fig.~\ref{Boxplotnevtncoal} shows that resampling among histories with the same number of lineages (that is by the number of coalescence events) is more efficient than the alternative condition. Indeed, for a given number of events, histories with fewer lineages (more coalescences) often correspond to a low current weight but might have a high final weight.

In Fig.~\ref{Boxplot_betafixe_alpha} and \ref{Boxplot_alphafixe_beta}, we propose to resample every  $k=8$, $6$, $2$ or $1$ coalescence events, which correspond respectively to resampling $12$, $16$, $49$ or $99$ times during Algorithm~3 since the simulated data set is composed of $100$ genes. We observe that the variation in likelihood estimation decreases with the frequency of the checkpoints, indicating that we should propose a resampling step as often as possible. Fig.~\ref{0.4_0_0.25_400_RMSEfctbeta_decalagencoal} and \ref{0.4_0_0.25_400_RMSEfctbeta_decalagealpha} leads to the same conclusion regarding the best calibration of $k$: MSE ratio comparing SISR with $k=1$ to the SIS estimate (represented by a $\ast$) are often the lowest MSE ratios when compared to other MSE ratios with the same tuning parameters $\alpha$ and $\beta$. Both Fig.~\ref{0.4_0_0.25_400_RMSEfctbeta_decalagencoal} and \ref{0.4_0_0.25_400_RMSEfctbeta_decalagealpha} also indicate that the MSE ratio always decreases with the frequency of checkpoints when the resampling distribution depends on the pairwise composite likelihood ($\beta = 0.01$).

\paragraph{Calibration of the powers $\alpha$ and $\beta$ of the resampling distribution}

The efficiency of the resampling distribution of Eq.~\eqref{NewWeigths} greatly differs between different values of the resampling parameters $\alpha$ and $\beta$. For this reason, they must be adequately chosen to balance the effect of the information provided by the SIS weights and by the composite likelihood.
First of all, we notice that the MSE Ratio is substantially below one for all sixty $(\alpha,\beta)$ couples considered (see Fig.~\ref{0.4_0_0.25_400_RMSEfctbeta_decalagencoal} and Fig.~\ref{0.4_0_0.25_400_RMSEfctbeta_decalagealpha}) meaning again that resampling improves the inference of the likelihood in terms of MSE, even when the resampling parameters are not much calibrated.
Then we observe that when the resampling distribution only depends on the SIS weights ($\beta=0$), the MSE Ratio globally decreases when $\alpha$ increases (see Fig.~\ref{0.4_0_0.25_400_RMSEfctbeta_decalagealpha}).
When resampling is performed after each coalescence, as recommended above, we find that any $\alpha$ between $0.5$ and $1$ is a reasonably good choice (Fig.~\ref{0.4_0_0.25_400_RMSEfctbeta_1coal}). This is further supported by  Fig.~\ref{Boxplot_ncoalfixe_beta} which represents the boxplots obtained with $k=1$, $\alpha=0.4$, $0.5$, $0.6$, $0.7$ and $1$ for three different values $\beta = 0$, $0.001$ and $0.01$. Indeed the variation in likelihood estimates corresponding to $\alpha=0.4$ is larger than for other values of $\alpha$ whereas the variation in likelihood estimates corresponding to $\alpha=0.7$ is slightly reduced compared to other $\alpha$ values. In the following we therefore set $\alpha=0.7$, which is higher than the arbitrary calibration $\alpha=0.5$ proposed by \citet{liu2001theoretical}.
Likewise with $k=1$, Fig.~\ref{0.4_0_0.25_400_RMSEfctbeta_1coal} represents the MSE ratios obtained with $\beta=0$ and $0.01$, and different values of $\alpha$. First we note that the ten values of MSE ratios are between $0.05$ and $0.20$, which means that the resampling reduces the MSE by a factor $5$ to $20$. Second the horizontal blue line indicates that the lowest MSE ratio obtained with $\beta=0$ is higher than the MSE ratio obtained with five out of seven values of $\alpha$ when $\beta=0.01$, that is when using the composite likelihood.
Moreover, this choice of $\beta$ is supported by Fig.~\ref{Boxplot_alphafixe_ncoal} which sets $\alpha=0.7$ and shows that the span of the distribution of the likelihood estimate is smaller when $\beta = 0.01$ (and $k=1$) than with other values of $\beta$.

\bigskip

In conclusion we can choose $k=1$, $\alpha = 0.7$ and $\beta=0.01$ although other choices of tuning parameters of the SISR also lead to an improvement of the likelihood estimation. The MSE ratio for the likelihood estimates using SIS vs. SISR is lower than $60\%$ for most choices of the resampling parameters, and below $10\%$ when the procedure proposes to resample after every coalescence event and when using the PCL, which allows a strong improvement of the likelihood estimation in a parameter point.

\section{Improvements in the likelihood based inference of demographic parameters}
 \label{sec:improvments}
The final aim of inference from a data set is not to compute the likelihood at a given point in parameter space. Rather, it is to provide an estimate of the parameter \(\phi=(\theta, D, \theta_\text{anc})\) and confidence intervals (CI) around each coordinates of $\phi$, i.e, marginal CIs. Thus, we will quantify the impact of our procedures on the performance of such inferences.

\subsection{Inference method}
\label{sub:inference}
In the numerical results below, we conduct a maximum likelihood analysis of the data. An estimate of the parameter $\phi$ is given by the maximum likelihood estimate (MLE) $\widehat{\phi}(\mathbf{x})$ of a multilocus data set $\mathbf{x}$. 
The multilocus likelihood is the product of the likelihoods for each locus $x_\ell$, $\ell=1, \ldots, d$:
\[
L_d(\phi, \mathbf{x}) = \prod_{\ell = 1}^d L(\phi, x_\ell),
\]
where $d$ is the number of loci in the sample, and each term of the product can be estimated with a SIS or a SISR algorithm.
The biological assumption that permits the above writing is that the loci are distant enough in the genome to have independent past histories.

Then we derived marginal CIs on each coordinate, from likelihood-ratio test based on the profile likelihood \citep[see][for details]{davison2003statistical}. To obtain a numerical value of the marginal CIs, we rely on asymptotic theory which states that under $H_0$ profile log likelihood-ratio is approximately $\chi^2$-distributed (e.g., \citealp{Severini2000}), here with a degree of freedom equal to $1$ when the size of the data set is large, i.e., when there is enough information on the parameter in the data.

\subsection{Inference algorithm and its evaluation} \label{InferenceAlgorithmAndItsEvaluation}

As in \cite{rousset2007likelihood}, \cite{rousset2012likelihood} and \cite{leblois2014maximum} we can conduct the inference process as follow. We first define a set of parameter points through a stratified random sample within a range of the parameter space provided by the user. Then, at each parameter point, the multilocus likelihood is the product of the likelihoods for each locus, which are estimated through the SIS or SISR algorithm.
The likelihoods inferred at the different parameter points are then smoothed by a Kriging scheme.
After a first analysis of the smoothed likelihood surface, the algorithm can be repeated a second time to increase the density of parameter points in the neighborhood of the first MLE.
The Kriging step removes part of the estimation error of likelihood in any given parameter point by assuming that the likelihood is a smooth function of the parameter. In this Section, we thus conducted numerical experiments to show that the gain of SISR over SIS in accuracy of likelihood estimates is retained through the Kriging step.

\bigskip

A convenient way to evaluate numerically the whole inference procedure, and in particular the coverage of the marginal CIs, is to check that distribution of $p$-value of the likelihood-ratio test of $H_0:\ \phi_1=\phi_1^\ast$ against $H_1:\ \phi_1\neq\phi_1^\ast$ is uniform on the interval $[0;\ 1]$ when the data set $\mathbf x$ is simulated from the model with $\phi_1=\phi_1^\ast$. To this end we represent the empirical cumulative distribution function (ECDF) of the $p$-value on many simulated data sets $\mathbf x$, which should be closed to the $1:1$ diagonal. Deviation from an uniform distribution can occur: either because the likelihood is poorly estimated or because the exact profile log likelihood ratios do not follow the asymptotic $\chi^2$ distribution. 
We perform a Kolmogorov-Smirnov test to assess the uniform distribution of the sample of $p$-values.

We have also computed other measures of the performance of the inference method of $\phi_1$ on many simulated data sets $\mathbf{x}_i$ from $\phi_1=\phi_1^\ast$, namely
\begin{itemize}
\item the mean relative bias of the MLE
, computed as
\[
\frac{(\text{observed bias on }\phi_1^\ast)}{\phi_1^\ast} = \frac{1}{\phi_1^\ast}\left(\operatorname{mean}_i\widehat{\phi}_1(\mathbf{x}_i)-\phi_1^\ast\right),
\]
\item and the relative root mean square error (relative RMSE) of the MLE
, computed as
\[
\sqrt{\frac{\text{MSE on }\phi_1^\ast}{\phi_1^{\ast 2}}} = \frac{1}{\phi_1}
\sqrt{\operatorname{mean}_i\left( \widehat{\phi}_1(\mathbf{x}_i)-\phi_1^\ast\right)^2}.
\]
\end{itemize}

\subsection{Numerical experiment cases and previous results}

We performed numerical experiments on four different demographic scenarios, all modeled according to the exponential contraction of Section~\ref{UnderlyingDemographicModel}, which are as follows.
\begin{enumerate}[label=\emph{(\roman*)}]
\item $\phi=(\theta, D,\theta_\text{anc})=(0.4,1.25,40)$, which is the baseline scenario of \citet{leblois2014maximum}: the population size has undergone a contraction of strength $N_\text{ratio}=0.01$ and $D=T/2N=1.25$, where $N$ is the size of the population at time $t=0$,
\item $\phi=(\theta, D,\theta_\text{anc})=(0.4,1.25,400)$ which differs from the baseline scenario in the strength $N_\text{ratio}=0.001$ of the contraction,
\item $\phi=(\theta, D,\theta_\text{anc})=(0.4,0.25,40)$ which differs from the baseline scenario in the speed of the contraction, since it corresponds to a contraction of strength $N_\text{ratio}=0.01$ but $D=T/2N=0.25$ which is five times smaller than in the baseline scenario,
\item $\phi=(\theta, D,\theta_\text{anc})=(0.4,0.25,400)$ which represents the strongest and recentest contraction of the four scenarios. Numerical results on the estimation of the likelihood in this last demographic scenario have already been presented in Section~\ref{Results} above.
\end{enumerate}

\bigskip

The baseline scenario is a case where the inference procedure performs well when the likelihood of each locus at each point of the parameter space is estimated with the sequential importance sampling on $n_H=2,000$ histories, but $n_H=100$ gives satisfactory results. 
A more careful look at the results shows that the profile likelihoods exhibit clear peaks around the MLE for all parameters. This is a first evidence that the data contain information regarding all parameters. And, indeed, the ECDF of the p-values are almost aligned on the $1:1$ diagonal (Fig.~\ref{0.4_0_1.25_40_50A_100A}).

\subsection{Numerical results}
To analyze the performances of the inference procedure we have simulated $500$ data sets for each scenario, composed of $10$ independent loci and $100$ genes per locus.

\paragraph{Gain in accuracy of the MLE}
The numerical results regarding the relative bias and relative RMSE are given in Table~\ref{tab:bias}. We find that the resampling technique allows overall to reduce the relative RMSE and the relative bias, mostly for the parameter $\theta$.

\begin{enumerate}[label=\emph{(\roman*)}]
\item For the baseline scenario $ (\theta=0.4, D=1.25, \theta_\text{anc}=40) $, the resampling procedure allows to reduce the relative bias on $\theta$ by $30\%$, the relative bias on $D$ by $20\%$ and the relative bias on $\theta_\text{anc}$ by $80\%$ and also the relative RMSE on $\theta_\text{anc}$ by $10\%$, other values being similar.

\item In the situation $ (\theta=0.4, D=1.25, \theta_\text{anc}=400)$ of a stronger contraction but not too recent, the resampling procedure allows to reduce the relative bias on $\theta$ by $40\%$,  and also the relative RMSE on $\theta$ by $35\%$, other values being similar.

\item Concerning the scenario $ (\theta=0.4, D=0.25, \theta_\text{anc}=40) $ of a more recent contraction than the baseline situation but with the same strength, the resampling procedure allows to reduce the relative bias on $\theta$ by $79\%$ and the relative RMSE on $\theta$ by $10\%$. However, the relative bias and relative RMSE on $D$ increase by a factor of $2.5$ and $1.8$ respectively.

\item We observe the same trend with the situation $ (\theta=0.4, D=0.25, \theta_\text{anc}=400) $ of a more recent and stronger contraction than the baseline situation. Indeed, when using resampling, the relative bias and relative RMSE on $\theta$ decrease by a factor $4.5$ and $3$ respectively while the relative bias and relative RMSE on the parameter $D$ increase by a factor $2.5$ approximately.
\end{enumerate}

\begin{table}[bth]
  \centering
  \caption{\bfseries Accuracy of the MLE with SIS and SISR algorithms}
  \label{tab:bias}
  \begin{tabularx}{\textwidth}%
{>{\bfseries}r @{\extracolsep{\fill}} r r @{\extracolsep{\fill}} r r @{\extracolsep{\fill}} r r @{\extracolsep{\fill}} r r @{\extracolsep{\fill}} r}
    \toprule
    \multicolumn{2}{r}{ \bfseries $(\theta, D,\theta_{anc})=$} &
    \multicolumn{2}{c}{$(0.4,1.25,40)$}&\multicolumn{2}{c}{$(0.4,1.25,400)$}  &
    \multicolumn{2}{c}{$(0.4,0.25,40)$} &\multicolumn{2}{c}{$(0.4,0.25,400)$}
    \\
    \midrule
    \multicolumn{2}{r}{\bfseries Algorithm} & SIS & SISR&SIS & SISR& SIS & SISR & SIS & SISR
    \\
    \multicolumn{2}{r}{\bfseries with $n_H=$} & \multicolumn{2}{c}{$100$} & \multicolumn{2}{c}{$100$} &
    \multicolumn{2}{c}{$2,000$} & \multicolumn{2}{c}{$2,000$}
    \\
    \midrule
    rel. bias& \(\theta\) & ${0.17}$  & ${0.12}$& $ {0.56} $ & ${0.34} $ & $ {1.07}$ & $ {0.23} $ &$ {4.9}$ & $ {1.1} $  \\
&D & $0.063$  & $0.051$& $ -0.02$ & $-0.017$& $ 0.23$ & $0.571$ & $ -0.06$ & $0.15$   \\
& \( \theta_{anc} \)& $0.076$  & $0.016$& $ 0.048$ & $-0.042$& $ 0.032$ & $0.023$ & $ 0.044$ & $-0.053$  \\
\midrule
rel. RMSE & \(\theta\) & ${0.47}$  & ${0.46}$& ${0.71} $ & ${0.53}$ & ${2} $ & ${1.8}$ & ${5.2} $ & ${1.7}$   \\
&D & $0.28$  & $0.28$& $ 0.14$ & $0.14$ & $ 0.45$ & $0.8$ & $ 0.14$ & $0.37$   \\
& \( \theta_{anc} \)& $0.53$  & $0.48$& $ 0.37$ & $0.38$ & $ 0.29$ & $0.27$ & $ 0.25$ & $0.23$  \\
\bottomrule
  \end{tabularx}
\end{table}

\paragraph{Reducing the size of the collection of histories thanks to the resampling}
We compare here the inference method relying on estimates of the likelihood either based on SIS or on SISR with a smaller collection of histories per locus and per point of the parameter space.

We first compare both methods on the baseline situation (\( \theta= 0.4, D =1.25, \theta_\text{anc} = 40 \)) of a relatively weak and not too recent contraction. In this case, the SIS algorithm performs well due to the large amount of information in the genetic data (\citealp{leblois2014maximum}).

We find out that SISR with \(50\) ancestral histories produces comparable results to SIS with \(100\) ancestral histories.  With the same number of ancestral histories, the relative bias and relative RMSE are lower with SISR than with SIS as explained above, and the ECDF of the $p$-values are closer to the diagonal as shown in Fig.~\ref{0.4_0_1.25_40_50A_100A}.  We conclude that resampling improves the parameter estimation in a situation where the SIS performs well, dividing by $2$ the required number of ancestral histories.

We then compare both methods on a more difficult situation (\( \theta= 0.4, D =1.25, \theta_\text{anc} = 400 \)) of a stronger contraction.
In this situation, the SIS procedure performs less well but we obtain satisfying results with \(n_H=2000\) sampled histories (\citealp{leblois2014maximum}), which leads to reasonable computation time. Indeed, for a single data set with one hundred gene copies and ten loci, analyses are carried out in $8$ hours in C++ process time on average. Here we decrease the number of ancestral histories in both SIS and SISR, in order to show that the SISR performs well when the SIS does not.
For both procedures we explored \(n_H=100\), \(200\) and \(400\) ancestral histories per parameter point. With $400$ ancestral histories, the analyses are carried out in  $1.7$ hours in C++ process time on average, for a single data set with one hundred gene copies and ten loci.

Fig.~\ref{0.4_0_1.25_400} shows that with the same number of explored ancestral histories, we obtain lower relative bias and relative RMSE and also better ECDF, which means better coverage properties of the CIs.
By comparing Fig.~\ref{0.4_0_1.25_400_200A_SIS} with \ref{0.4_0_1.25_400_100A_SISR}, we also find that the SISR procedure provides comparable results (relative bias, relative RMSE and ECDF) to the SIS procedure with half of the number of explored ancestral histories, as in the baseline scenario. The comparison of Fig.~\ref{0.4_0_1.25_400_400A_SIS} with Fig.~\ref{0.4_0_1.25_400_200A_SISR} shows the same results.

\paragraph{Improvements in more difficult situations}
We compare both SIS and SISR procedure on two more difficult situations $ (\theta=0.4, D=0.25, \theta_\text{anc}=40) $ and $ (\theta=0.4, D=0.25, \theta_\text{anc}=400)$ of very recent contraction, very strong for the second case. Thus in these two cases, the proposal distribution is inefficient. In the extreme case of $ (\theta=0.4, D=0.25, \theta_\text{anc}=400)$, the SIS algorithm does not provide satisfactory CI coverage properties, even with $200,000$ sampled histories.

We first consider the situation (\( \theta= 0.4, D =0.25, \theta_\text{anc} = 40 \)) in which the magnitude of the contraction is the same as the baseline case but here it occurred much more recently.
In this situation, the SIS procedure performs less well and we do not obtain satisfactory results with \(2000\) sampled histories (\citealp{leblois2014maximum}).  With the same number of ancestral explored ancestral histories, applying the resampling technique, we obtain lower relative bias and relative RMSE for \(\theta\) but higher for \(D\). We also obtain less good CI coverage properties for \(\theta\) and \(D\) (see Fig.~\ref{0.4_0_0.25_40_2000A}). 
We conclude that the method detects the contraction, estimates the date and the ancestral size of the population but does not find enough information in the data about the current size of the population. Indeed, the likelihood surface is quite flat. 

We then consider the extreme situation with a very recent and stronger past contraction (\( \theta= 0.4, D =0.25, \theta_\text{anc} = 400 \)). Increasing greatly the number of ancestral histories sampled per parameter point up to \(200,000\) and consequently the computation times by a factor $100$, decreases relative bias and relative RMSE on \(\theta\) but does not provide satisfactory CI coverage (\citealp{leblois2014maximum}). The resampling technique allows us to decrease relative bias and relative RMSE of \(\theta\) for a fixed number of ancestral histories sampled and also to provide better CI coverage (see Fig.~\ref{0.4_0_0.25_400}). Comparing Fig.~\ref{0.4_0_0.25_400_200000A_SIS} and \ref{0.4_0_0.25_400_2000A_SISR}, we obtain approximately the same relative bias and relative RMSE on \(\theta\) and close CI coverage performance with SISR with \(100\) times less histories sampled than with SIS.
 In this scenario as in the previous one there is not enough information in the data and hence the likelihood surface is flat. We also face an issue due to the strength of the contraction and the inefficiency of the proposal distribution in this situation, in addition to the previous difficulty. The resampling technique allows a global gain, since it corrects the inefficiency of the proposal distribution in disequilibrium situations. However, in this situation as in the previous one, the resampling technique does not bring an improvement to the lack of information in the data about the parameter $\theta$ when the contraction is too recent.

\paragraph{About flat likelihood surfaces}
The inference method of Section~\ref{sub:inference} suffers from two major defects when the likelihood surface is flat. First any error regarding the likelihood estimate at some point of the parameter space can lead to artificial local maxima or minima in the smoothed surface, depending on how the Kriging method performs. Second even if we recover the true flat likelihood surface, the distribution of the likelihood-ratio statistic is not correctly approximated by the $\chi^2$ distribution because of the lack of information on the parameter in the data. Both defects can be seen in our numerical studies. First the median of the SIS likelihood estimates is much lower than the reference value in
Fig.~\ref{Boxplotnevtncoal} to \ref{Boxplot_alphafixe_ncoal}. It indicates that the SIS likelihood estimates are very often below the true value of the likelihood which can lead to artificial local minimum in the smoothed likelihood surface obtained with the Kriging algorithm. When comparing Fig.~\ref{0.4_0_0.25_40_2000A_profiles} (a) and (b), the SIS likelihood estimates introduce a local minimum at $\theta=2N\mu\approx 0.001$ while the more reliable SISR likelihood estimates manage to recover the flat likelihood surface. Likewise on $D=T/2N$, the SIS likelihood estimates introduce a local minimum around $D\approx 0.4$, see Fig.~\ref{0.4_0_0.25_40_2000A_profiles} (c), while the SISR likelihood estimates recover the flat likelihood surface, see Fig.~\ref{0.4_0_0.25_40_2000A_profiles} (d).
The second defect can be seen on the p-value's ECDF in Fig.~\ref{0.4_0_0.25_400}. Even with the SISR likelihood estimates which manage to recover the flat likelihood surface, the Figure exhibits a departure from the uniform distribution, meaning that the $\chi^2$ approximation is not accurate.
Flat likelihood surface as in Fig.~\ref{0.4_0_0.25_40_2000A_profiles} and very large (even if untrustworthy) CIs, for instance the $95\%$ CI with SISR likelihood estimates is approximately $[0.00016;  2.5]$ on $\theta = 2 N \mu$  and $[0.14;0.72]$ on $D = T / 2 N$, should act as warnings that the data do not carry much information about the parameters of interest.

\section{Cynopterus sphinx data set}

We applied the inference method of \ref{sub:inference} on the fruit bat \textit{Cynopterus sphinx} data set presented in \cite{storz2002testing},
consisting in allelic frequencies computed from a sample of $246$ individuals, hence $492$ genes, genotyped at $8$ microsatellite loci.
Using a coalescent-based Monte Carlo Markov chain (MCMC) algorithm, \cite{storz2002testing} found a strong evidence for a pronounced population contraction
with (i) a posterior mode for $N_{ratio}$, the strength of the contraction, around 0.03 (but with large $90\% $ highest probability density (HPD) interval of $[0.0005 - 0.5]$);
and (ii) a mode for the scaled time of occurrence, $D=T/2N$, around $0.35$ with $90\% $ HPD interval of $[0.015 - 0.9]$.

Here, we compare their results with those obtained with our methods, with an emphasis on the gain in precision and in computation cost
using the SISR algorithm with different resampling probabilities vs. the SIS algorithm.
To assess the accuracy of the different algorithms on this real data for which we do not have concrete expectations, we rely on the following steps.

We know from the results presented in the previous sections that both SIS ans SISR algorithms may give inaccurate point estimate and CIs
for some parameters if the number of simulated histories is too low to infer the likelihood at each parameter point with enough precision.
For this reason, we first compared parameter estimates obtained by analyzing the data with the SIS algorithm using $n_H=1,000; 10,000; 100,000$ and $1,000,000$ sampled histories per point.
Moreover, to evaluate the IS variance of the log-likelihood estimate by the important sampling algorithm, we also considered the likelihood RMSE estimate from kriged duplicate points. It is computed from independent pairs of likelihood estimates at the same parameter point
for the points retained at the Kriging step.

\begin{table}[bth]
  \small{
  \centering
  \caption{\bfseries Evolution of the MLE with SIS when increasing the number of sampled histories per point, on the bat data set}
  \label{tab:CynopterusSphinxSIS}
  \begin{tabularx}{\textwidth}%
{>{\bfseries} l @{\extracolsep{\fill}} cccc @{\extracolsep{\fill}} c  @{\extracolsep{\fill}} c  @{\extracolsep{\fill}}  }

    \toprule
    \multicolumn{1}{l}{ \bfseries $n_H$} &
    \multicolumn{1}{c}{$ \hat{\theta}$}&
    \multicolumn{1}{c}{$ \hat{D} $}&
    \multicolumn{1}{c}{$\widehat{\theta_\text{anc}}$}&
    \multicolumn{1}{c}{$ \widehat{N_\text{ratio}} $}&\multicolumn{1}{c}{$\ln(\hat{L}(\hat{\phi}))$}
    &   \multicolumn{1}{c}{lik-RMSE}
     \\
    \midrule
    $1,000$& $1.1$  & $0.44$& $320$&$0.0034$&$-663.74$& $9.2$  \\
    &$[0.76 - 1.5]$ & $[0.35 - 0.55]$ & $ [180 - 560]$ &$ [0.0018 - 0.0062]$ &   \vspace{5pt} \\
     $10,000$& $0.80$  & $0.45$& $340$&$0.0024 $&$-639.08$ & $8.3$   \\
     &$[0.53 - 1.16]$ & $[0.37 - 0.55]$ & $ [205 - 553]$ &$ [0.0013 - 0.0045]$ & \vspace{5pt} \\
     $100,000$& $0.62$  & $0.48$& $460$&$0.0014$&$-620.68$ & $5.7$  \\
     &$[0.38 - 0.94]$ & $[0.40 - 0.57]$ & $ [254 - 492]$ &$ [0.00078 - 0.0027]$ &   \vspace{5pt}\\
     $1,000,000$& $0.42$  & $0.49$& $380$&$0.0011$&$-606.59$ & $4.4 $  \\
     &$[0.24 -0.68]$ & $[0.41 - 0.59]$ & $ [240 - 589]$ &$ [0.00056 - 0.0022]$ &  \\

  \midrule
  \end{tabularx}
  }
\end{table}

Our results, presented in Table~\ref{tab:CynopterusSphinxSIS}, show a strong decrease in the MLEs of $\theta$ et $N_\text{ratio}$, from $1.1$ to $0.42$ and $0.0034$ to $0.00111$ respectively , with an associated shift and narrowing of CIs.
An opposite and weaker trend of increase is also observed for $\hat{D}$, from $0.44$ to $0.49$, and for $\widehat{\theta_\text{anc}}$, from $320$ to $380$, with a slight shift of CIs for both parameters.
As expected, increasing the number of simulated histories leads to a decrease in the variance of the likelihood estimate at each parameter point, but more interestingly,
it also leads to an important increase of $\ln(\hat{L}(\hat{\phi}))$ the log-likelihood value at the MLE, from $-664$ to $-607$.
This first step implied a very high computational burden but allowed us to see a clear trend towards an improvement of the SIS inference for each analyses with more histories.
We probably did not reach a plateau in the precision even when considering $n_H=1,000,000$ histories (which took approximately $15000$ hours in C++ process time),
but considering more histories is merely unfeasible. We thus consider the results obtained with $1,000,000$ histories as the best result we can get with the SIS algorithm,
irrespectively of computation costs.

In a second step, we compared estimations obtained with the SIS and various versions of SISR algorithms by playing with different resampling strategies,
for a fixed and reasonable computation cost.
To this end, we fixed the number of sampled histories per parameter point $n_H$ to 1,000, leading to an average of $15$ hours in C++ process time for each analysis.
Table~\ref{tab:CynopterusSphinx} presents, for each analysis, the point estimates and CIs for each parameter of interest, including the $N_\text{ratio}$,
as well as the log-likelihood value estimated at the MLE $\ln(\hat{L}(\hat{\phi}))$ and the RMSE estimate for the log-likelihood at each parameter point (i.e. the IS log-likelihood variance).

A first global conclusion is that all SISR analyses with different resampling probability distributions lead to very similar results both in terms of parameter estimates and maximum log-likelihood value over the surface. Moreover, results from SISR analyses differ from those obtained with the SIS algorithm with the same number of sampled histories, but are closer to the best results we get with the SIS algorithm when considering 1,000,000 histories. Our results thus show that SISR analyses on this data set lead to more accurate inferences than SIS for similar computation costs.
According to Table~\ref{tab:CynopterusSphinx}, the most significant improvement due to resampling concerns the parameter $\theta$, which is the most poorly estimated by the SIS algorithm with $n_H=1,000$.
To a lesser extent, resampling also improves inference of the parameter $D$ compared to the SIS algorithm, but conclusions about the parameter $\theta_\text{anc}$ are less obvious, even if almost all analyses lead to values closer to the best estimate we have than the one obtained without resampling.

Tackling precisely the effects of the different resampling strategies on this real data set analysis is more difficult as we do not have in hand the true values of the different parameters
and also because our results show subtle differences among parameters for different resampling strategies. Thus, we can only globally conclude that resampling, especially when using PCL, always lowers the likelihood estimation variance and increases the maximum likelihood value on this data set compared to the analysis without resampling. It should also be noted that resampling according to a distribution which does not depend on the IS weights but only on the PCL also improves the inference precision compared to SIS in terms of likelihood estimation variance, maximum likelihood value and parameter estimates (results not shown).

Overall, our results on this fruit bat data set analyses qualitatively agree with those of \cite{storz2002testing} showing a strong evidence for a past population contraction.
However, comparison between their results and our own inferences shows important quantitative differences in point estimates and associated CIs.
First, we found an almost hundred times stronger contraction (e.g. $\widehat{N_\text{ratio}}$ around 0.0005 vs. 0.03) that occurred slightly further in the past with our $\hat{D}$ being about 0.55 vs. 0.35.
Moreover, our analyses show greater precision for those parameter estimates with CIs that are much narrower than their credibility intervals.
Second, despite \cite{storz2002testing} do not present their results in terms of $\theta$ and $\theta_\text{anc}$, we found that our scaled estimates of current population size (i.e. $\theta$)
are congruent with their unscaled estimates of current population sizes and mutation rates.
The large differences in the amplitude of the past contraction inferred by \cite{storz2002testing} thus comes mostly from differences in past population size estimates.
Our SIS and SISR analyses inferred ancestral population sizes that are hundred times larger than their estimates.
This latter result do not seem to be due to a potential bias in our methods as ancestral population size estimates do not decrease with increasing precision (low RMSE) of estimates of likelihood (the opposite may actually occur).

\begin{table}[bth]
 \footnotesize{
  \centering
  \caption{\bfseries Accuracy of the MLE with SIS and SISR algorithms on the  bat data set}
  \label{tab:CynopterusSphinx}
  \begin{tabularx}{\textwidth}%
{>{\bfseries} l  @{\extracolsep{\fill}} c c c c@{\extracolsep{\fill}} c @{\extracolsep{\fill}} c  }
    \toprule
    \multicolumn{1}{l}{ \bfseries $(\alpha,\beta,k)$} &
   \multicolumn{1}{c}{$ \hat{\theta}$}&
    \multicolumn{1}{c}{$ \hat{D} $}&
    \multicolumn{1}{c}{$\widehat{\theta_\text{anc}}$}&
    \multicolumn{1}{c}{$ \widehat{N_\text{ratio}} $}&\multicolumn{1}{c}{$\ln(\hat{L}(\hat{\phi}))$}  &   \multicolumn{1}{c}{lik-RMSE}
     \\
    \cmidrule{1-7}
   \multirow{2}{*}{\mbox{SIS}}  & $1.08$  & $0.44$& $320$&$0.0034$& $-663.74$ &  $ 9.2 $  \\
 &$[0.76 - 1.5]$ & $[0.35 - 0.55]$ & $ [180 - 560]$ &$ [0.0018 - 0.0062]$ &  &  \\

    \cmidrule{1-7}
     \multirow{2}{*}{\mbox{SISR}} &&&&&& \vspace{5pt}\\
   \multirow{2}{*}{\mbox{$(1,0,50)$}} &$0.20$  & $0.54$& $490$&$0.00041$&$-586.81$&  $ 6.7 $  \\
    &$[0.085 - 0.4]$ & $[0.46 - 0.62]$ & $ [270 - 490]$ &$ [0.00017 - 0.00098]$ &  & \vspace{5pt} \\

\multirow{2}{*}{\mbox{$(1,0.0005,50)$}} &$0.18$  & $0.53$& $380$&$0.00047$&$-587.38$&  $ 5.8 $ \\
 &$[0.073 - 0.37]$ & $[0.45 - 0.63]$ & $ [240 - 500]$ &$ [0.00018 - 0.0011]$ &  &  \vspace{5pt}\\

\multirow{2}{*}{\mbox{$(1,0.001,50)$}} &$0.21$  & $0.53$& $390$&$0.00053$& $-586.78$&  $5.3 $\\
 &$[0.091 - 0.39]$ & $[0.45 - 0.62]$ & $ [250 - 490]$ &$ [0.00022 - 0.0011]$ &  &  \\
 \cmidrule{1-7}

\multirow{2}{*}{\mbox{$(0.7,0,50)$ }}&$0.19$  & $0.53$& $410$&$0.00046 $&$-588.25$&  $6.5 $ \\
 &$[0.075 - 0.40]$ & $[0.45 - 0.63]$ & $ [230 - 510]$ &$ [0.00016 - 0.0012]$ &  &  \vspace{5pt}\\

\multirow{2}{*}{\mbox{$(0.7,0.0005,50)$}} &$0.20$  & $0.54$& $450$&$0.00046$&$-588.11$&  $6.2$\\
 &$[0.084 - 0.43]$ & $[0.45 - 0.63]$ & $ [250 - 490]$ &$ [0.00018 - 0.0011]$ &  & \vspace{5pt} \\

\multirow{2}{*}{\mbox{$(0.7,0.001,50)$ }}&$0.22$  & $0.53$& $440$&$0.00050$&$-588.25$&  $8.0$ \\
 &$[0.093 - 0.43]$ & $[0.45 - 0.63]$ & $ [272 - 581]$ &$ [0.00019 - 0.0011]$ &  &  \\
 \cmidrule{1-7}

\multirow{2}{*}{\mbox{$(0.5,0,50)$}} &$0.21$  & $0.53$& $400$&$0.00053$&$-588.43$&  $ 6.4 $ \\
 &$[0.088 - 0.40]$ & $[0.44 - 0.61]$ & $ [240 - NA]$ &$ [0.00020 - 0.0012]$ &  &  \vspace{5pt}\\

\multirow{2}{*}{\mbox{$(0.5,0.0005,50)$}} &$0.17$  & $0.52$& $360$&$0.00048 $&$-589.07$&  $5.0$ \\
&$[0.066 - 0.38]$ & $[0.44 -- 0.61]$ & $ [230 - 490]$ &$ [0.00017 - 0.0012]$ &  &  \vspace{5pt}\\

\multirow{2}{*}{\mbox{$(0.5,0.001,50)$}}& $0.21$  & $0.54$& $490$&$0.00042$&$-588.99$& $ 7.4$ \\
 &$[0.091 - 0.40]$ & $[0.46 - 0.61]$ & $ [330 - 490]$ &$ [0.00018 - 0.00089]$ &  &  \\
 \cmidrule{1-7}
  \end{tabularx}
  }
\end{table}

\section{Conclusions}

In this study, we proposed to improve sequential importance sample algorithms for likelihood inference of demographic parameters
by adding a resampling procedure based on \citet{liu2001theoretical}. The new resampling probability distribution we considered depends on the SIS weights,
as proposed in \citet{liu2001theoretical} and \citet{liu2008monte}, but also on the pairwise composite likelihood of the sample,
providing additional information about the future trend of each sampled history. To evaluate the gain in efficiency due to this new strategy,
we focused on varying population size models for which no efficient proposal distribution is available.

We first showed using simulations that resampling allows to reduce the variance and the bias in the likelihood estimate at a parameter point. This step allowed us to show that resampling is more efficient when the checkpoints correspond to a number of coalescent events undergone by each history rather than considering coalescent and mutations events together, and that checkpoints must be frequent, more precisely after each coalescence event. This step also showed that the information provided by the composite likelihood allows a stronger decrease in variance of the likelihood estimates than the information provided only by the SIS weights. 

Then we showed that the increased precision in estimation of likelihood also improved the likelihood-based inferences. According to numerical results based on the analysis of simulated data sets under various scenarios of past population size contraction with different strength and timing, we conclude that the resampling procedure helps to correct the inefficiency of the SIS proposal distribution. Under those time-inhomogeneous models, the stronger the population size contraction is, the less efficient the proposal distribution is because it is computed under equilibrium hypothesis. For a similar precision of inference, SISR provides al least a two-fold, and sometimes a much higher gain in computational efficiency for these models. However, when the contraction is very recent, the genetic data does not contain much information about the parameter $\theta$, resulting in flat likelihood surfaces. 
Consequently, our simulation results show that (i) when the contraction of the population size is not very recent, the resampling procedure divides the computation cost by at least a factor $2$; (ii) when the contraction is not very strong but very recent, it is not obvious how the resampling improves the inference because we face an important lack of information for the $\theta$ parameter; and finally (iii) in extreme cases of both very recent and very strong contractions, the resampling procedure divides the computation cost by a factor $100$, partially correcting the inefficiency of the SIS proposal distribution but not the lack of information in the data, and thus leads to potential bias in $\theta$ estimates and associated incorrect CI coverage properties.

Finally, we analyzed a fruit bat microsatellite data set previously analyzed by \cite{storz2002testing} using coalescent-based MCMC algorithms and found that the SISR method allows to infer the past demographic history of this population with computation times reduced by at least a factor one hundred compared to SIS. Our results also shows a much more pronounced past contraction of population size, with larger ancestral ancestral size, than previously found.

\section{Acknowledgments}

The first author were financially supported by the Labex Numev (Solutions Num\'eriques, Mat\'erielles et
Mod\'elisation pour l'Environnement et le Vivant, ANR-10-LABX-20) and by the Labex CeMEB (Centre M\'editerrann\'een de l'Environnement et de la Biodiversit\'e), Projet Investissement d'Avenir 2011-2019. Part of this work was carried out by using the resources of the INRA MIGALE (http://migale.jouy.inra.fr) and GENOTOUL (Toulouse Midi-Pyr\'en\'ees) bioinformatics platforms, the computing grids of the CBGP lab and the Montpellier Bioinformatics Biodiversity platform.

\bibliographystyle{plainnat}
\bibliography{mabiblio}

\appendix
\section{Pairwise Composite Likelihood}
\label{AppendixPCL}

In this work, we propose a new resampling probability distribution $\mathbf{v}$ on the collection of histories. 
To achieve efficiency we should pick a distribution $\mathbf{v}$ that reflects the future trend of the partial histories, more precisely that sets relatively high probabilities on the partial histories that will be the most likely at the end of the simulation.\\
Since the information is mainly in the dependency between individuals of the sample, we propose to substitute the Pairwise Composite Likelihood (PCL), noted $L_2(\mathbf{h})$ , for the likelihood of the end of the history $\mathbf{h}$. Indeed, it is the product of the likelihoods of each pair of individuals remaining in the sample :
\begin{equation*}
L_2(\mathbf{h};\phi) = \prod _ {A \in E: \, \mathbf{h}(A) \geq 2} L_2((A,A);\phi)^{\mathbf{h}(A)(\mathbf{h}(A)-1)/2} \times \prod _ {A < B \in E: \, \mathbf{h}(A)\geq 1, \, \mathbf{h}(B) \geq 1}L_2((A,B);\phi)^{\mathbf{h}(A)\mathbf{h}(B)},
\label{PCL}
\end{equation*}
where  $L_2(x,y|\phi)$ is the likelihood of a sample of two alleles in the stepwise mutational model, known explicitly when the scaled population size is constant equal to $\theta_\text{anc}$:
\begin{equation*}
\forall x,y \in E, \quad L_2((x,y);\phi) = \frac{1}{\sqrt{1+2\theta_\text{anc}}}\rho(\theta_\text{anc})^{|y-x|} \qquad \text{ where: } \qquad
\rho(\theta_\text{anc}) = \frac{\theta_\text{anc}}{  1 +\theta_\text{anc} + \sqrt{1 + 2\theta_\text{anc}}}.
\label{43}
\end{equation*}

The PCL does not approximate the likelihood but it behaves the same way. In particular, it reflects the behavior of the remaining lineages in the history $\mathbf{h}$ when the population size is equal to the ancestral size $N_\text{anc}$. Thus we chose to compute the PCL as an expression depending on $\phi$ only through $\theta_\text{anc}$. A judicious idea could be to derive an expression of the PCL with $\theta(t)$ depending on $t$. However, we are concerned that its calculation would become too complex with respect to the expected gain.


\section{Simulating holding times of $\mathbf{H}$ and $\widetilde{\mathbf{H}}$}
\label{AppendixGillespie}

The importance sampling scheme assumes that we are always able to draw simulations from the holding time
distribution. When the process $\widetilde{\mathbf{H}}$ has just jumped to state $\mathbf{h}\in \mathbb{N}^E$ at time
$t$, the distribution of the  holding time $\Delta$ has the following density
\[
\mathbb{P}\big(\Delta\in(\delta;\delta+d\delta)\big|T_i=t,\mathbf{X}_i=\mathbf{h}\big) =
\lambda_{t+\delta}(\mathbf{h}) \exp\left( -\int_0^\delta \lambda_{t+u}(\mathbf{h})du\right)d\delta
\]
and the following probability distribution function
\[
\mathbb{P}\big(\Delta\le\delta\big|T_i=t,\mathbf{X}_i=\mathbf{h}\big) =
1 - \exp\left( -\int_0^\delta \lambda_{t+u}(\mathbf{h})du\right).
\]
The intensity $\lambda_t(\mathbf{h})$ is defined by Eq.~\eqref{eq:lambda} and depends on the parametric model 
we set on the population size $N(t)$.
When the population size is constant, $\lambda_t(\mathbf{h})$ does not depend on $t$ and this distribution
boils down to the simple exponential distribution with rate $\lambda(\mathbf{h})$.
When we face an exponentially changing population size, the above probability distribution function can be
computed explicitly, see \citet{griffiths1994sampling}, and we can rely on the inverse of the probability
distribution function to draw simulations.

\begin{table*}[hbt]
\hrule
\vspace{1mm}
\textbf{Algorithm~4}: Simulation of the holding time
\vspace{1mm}
\hrule
\begin{enumerate}
\item Initialization: set $t' = t$
\item  \textbf{Repeat}
          \vspace*{-3mm}
\begin{enumerate}[label=2.(\alph*)]
        \item Compute some $M$ greater or equal to $\sup_{s\ge t'}\lambda_s(\mathbf{h})$
        \item  Draw $\delta_0$ from the exponential distribution with rate $M$
        and Set $t'=t'+\delta_0$
        \item  Draw $u$ uniformly in $[0,1)$
\end{enumerate}
        \hspace*{-.7cm} \textbf{until} $u \le \lambda_{t'}(\mathbf{h})\big/M$ 
       \item  Return $t'-t$
\end{enumerate}
\hrule
\end{table*}

In the general case where $N(t)$ can be any parametric function of $t$, Algorithm~4 might be interesting. It
follows the well known algorithm of \citet{gillespie1977exact} which aims at simulating pure jump Markov
processes. There is many ways to show that Algorithm~4 is correct. The simplest one is to claim that the holding
time is exactly the first point after $t$ of a Poisson point process \citep{kingman1992poisson} with intensity
$\lambda_s(\mathbf{h})$ at any $s\ge t$. And such Poisson point processes can be simulated from a first
Poisson point process with intensity $M$, constant over time, bounding the intensity $\lambda_s(\mathbf{h})$
for any $s\ge t$. To this end, we simply takes any point $T'$ from the first process with probability
$\lambda_{T'}(\mathbf{h})/M$. Thus, Algorithm~4 simulates the first point $t'$ of Poisson point process with
intensity $M$ at step 2.(b) and rejects it until an event of probability $\lambda_{t'}(\mathbf{h})/M$ occurs,
namely that $u\le \lambda_{t'}(\mathbf{h})/M$.  Finally, note that we can easily bound from above the
intensity $\lambda_s(\mathbf{h})$ at any $s\ge t$ by bounding from below the population size after time $t$, see
Eq.~\eqref{eq:lambda}.

\section{One-parameter profile likelihood ratios}
\label{ProfileLikelihood}
The likelihood function depends here on three parameter $( \theta, D, \theta_\text{anc})$. 
For each fixed value of the parameter $\theta$, we rewrite the likelihood estimate surface as $\hat{L}_{\theta}(D,\theta_\text{anc})=\hat{L}(\theta, D,\theta_\text{anc})$ that is $\theta$ is fixed and $D$ and $\theta_\text{anc}$ vary. 
We estimate $D$ and $\theta_\text{anc}$ by maximizing $\hat{L}_{\theta}(D,\theta_\text{anc})$ with respect to $D$ and $\theta_\text{anc}$, {\it i.e.}
\[(\hat{D}_\theta,\hat{\theta}_{\text{anc }\theta}) = \text{argmax}_{D, \theta_\text{anc}} \hat{L}_{\theta}(D,\theta_\text{anc}).\]
 As $\theta$ is unknown, we evaluate $(\hat{D}_\theta,\hat{\theta}_{\text{anc }\theta})$ for each $\theta$. 
 Then we estimate $\theta$ by maximizing $\hat{L}_{\theta}(\hat{D}_\theta,\hat{\theta}_{\text{anc},\theta})$ : 
 \[\hat{\theta} = \text{argmax}_{\theta} \hat{L}(\theta,\hat{D}_\theta,\hat{\theta}_{\text{anc},\theta}).\]
We have profiled out the parameters $(\hat{D},\hat{\theta}_\text{anc})$ and the  likelihood profile $\hat{L}_{\theta}(\hat{D}_\theta,\hat{\theta}_{\text{anc},\theta})$ is completely in terms of the parameter $\theta$. 
Then, we represent the likelihood profile function divided by its maximum value : \[\theta \mapsto \hat{L}_{\theta}(\hat{D}_\theta,\hat{\theta}_{\text{anc }\theta})/\hat{L}(\hat{\theta},\hat{D}_\theta,\hat{\theta}_{\text{anc},\theta}),\] 
We also represent in the same way the profile likelihood ratio for $D$ and the profile likelihood ratio for $\theta_\text{anc}$, estimated here with an IS algorithm.

\section{Figures}
\begin{figure}[!h]
  \centering
  \def\svgwidth{.8\textwidth}
  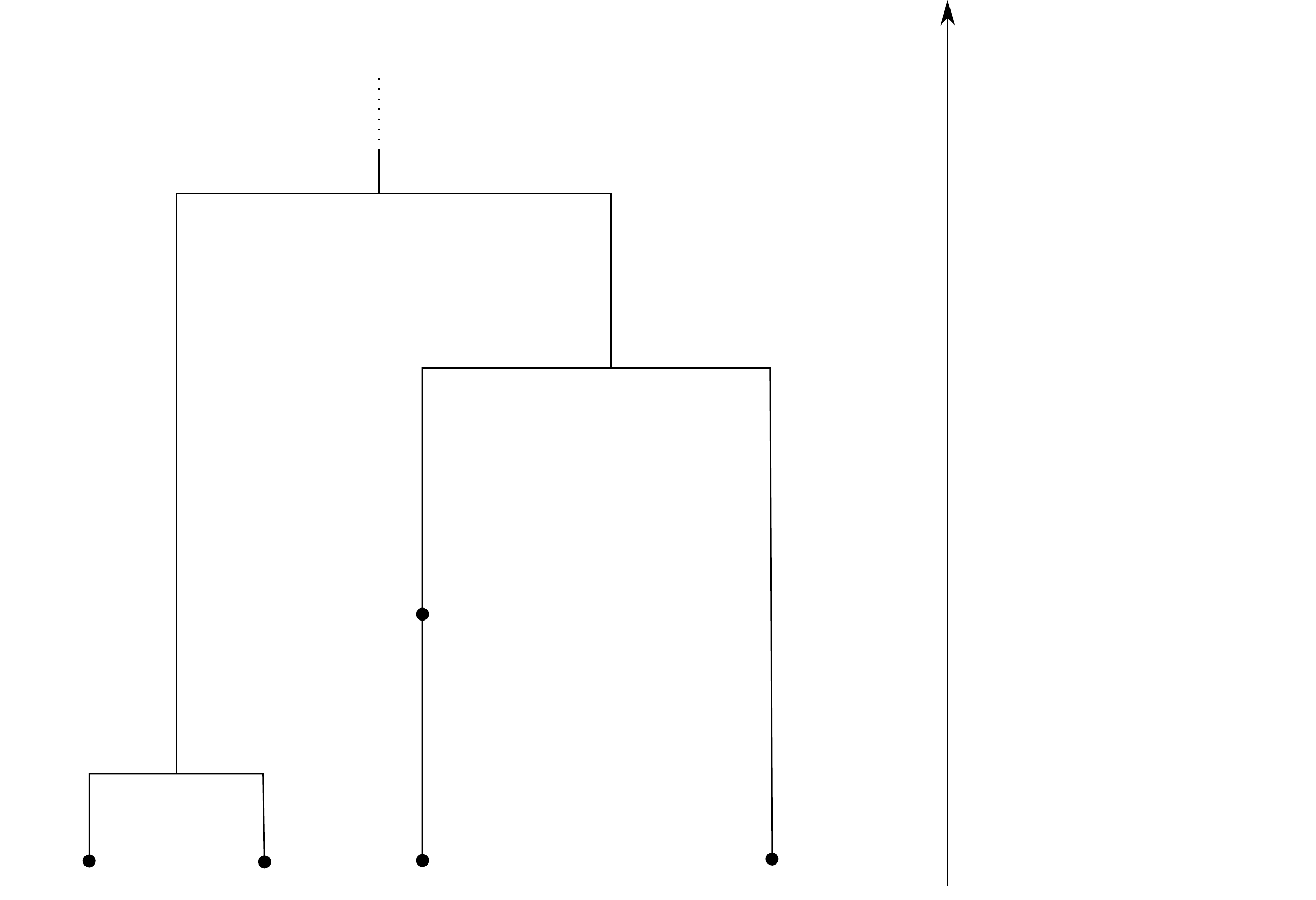
  \caption{An example of path of the process $\mathbf{H}$ from the MRCA leading to a sample of 4 genes at time $0$, when the set of gene types 
    $E=\{\text{a}, \text{b}, \text{c}\}$ is composed of three possible alleles.}
  \label{fig:coal}
\end{figure}


\begin{figure}[!h]
\begin{center}
\includegraphics[scale=0.85]{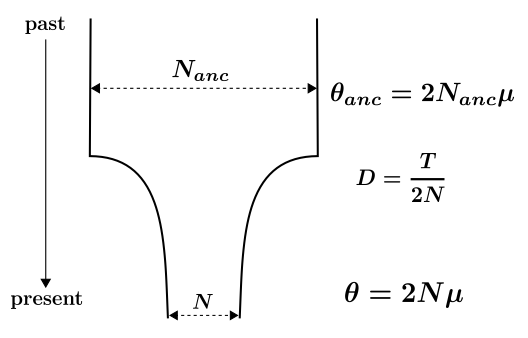} 
\caption{\textbf{Demographic model.}
 Representation of the exponentially contracting population size model used in the study. $N$ is the current population size, $N_\text{anc}$ is the ancestral population size (before the demographic change), $T$ is the time measured in generation since present, and $\theta$ is the mutation rate of the marker used. Those four parameters are the canonical parameters of the model. $\theta$, $D$, and $\theta_\text{anc}$ are the inferred scaled parameters (\cite{leblois2014maximum}).}
\label{ECP}
\end{center}
\end{figure}


\begin{figure}[h!]
\centering
\includegraphics[scale=0.63]{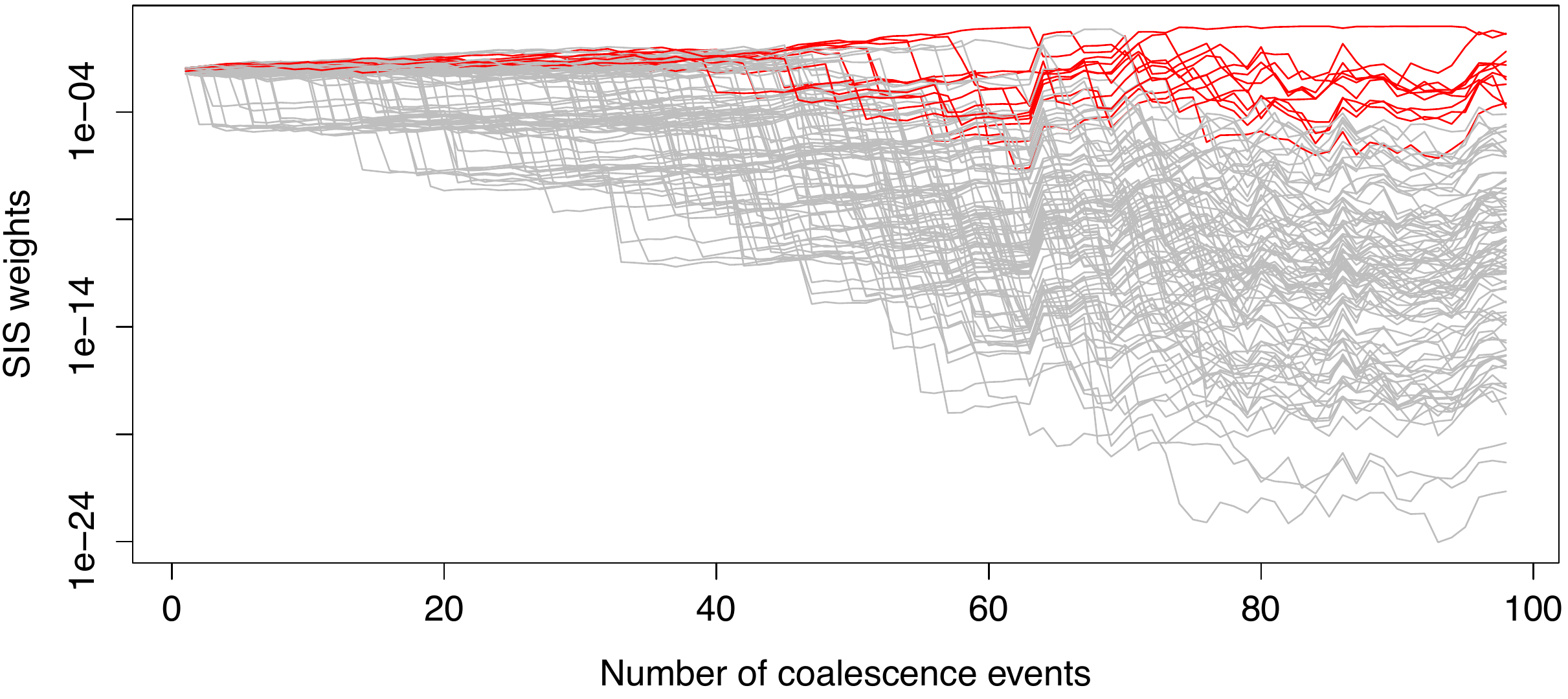}
\caption{\textbf{Evolution of the current partial SIS weights as a function of the number of coalescence events undergone by the sample}, for $100$ histories drawn according to the importance distribution.
At each step, the current weights are normalized by the sum of the current weights of all the histories and represented on a logarithmic scale. In red: the $10$ histories contributing the most to the overall estimate of the likelihood.}
\label{0.4_0_0.25_400_NormaliseWeightTrajectory}
\end{figure}

\begin{figure}[!h]
   \centering
\includegraphics[scale=0.7]{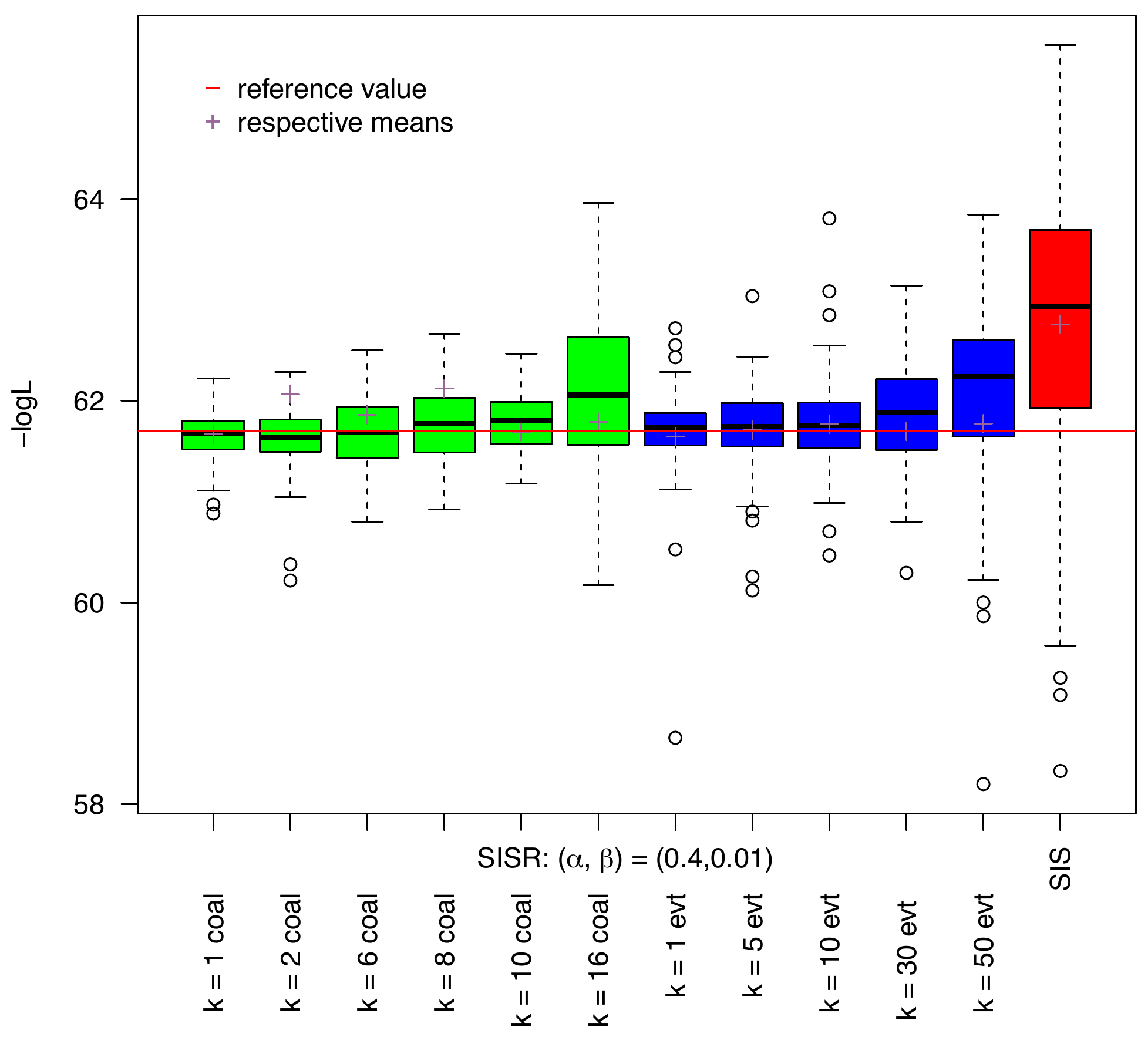}
\caption{\textbf{Boxplots of \(100\) estimates of the likelihood in a given parameter point with different inference algorithms.} The green boxes were obtained by resampling among histories with the same number of coalescences undergone while the blue boxes were obtained by resampling among histories with the same number of events undergone and the red box correspond to SIS estimates.}
\label{Boxplotnevtncoal}
\end{figure}

\begin{figure}[!h]
\centering
\subfloat[]{
\includegraphics[scale=0.4]{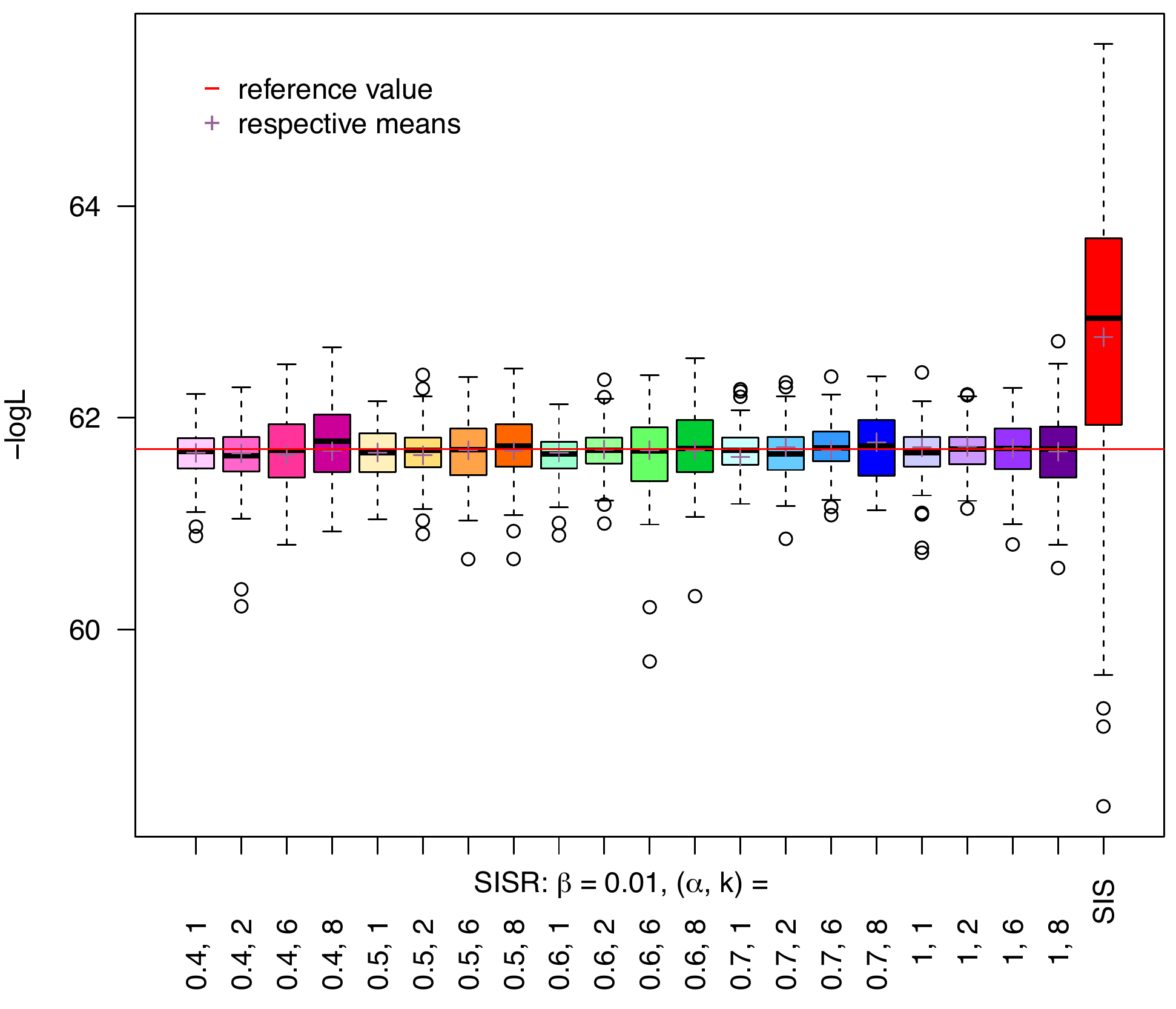}\label{Boxplot_betafixe_alpha}}
\subfloat[]{
 \includegraphics[scale=0.385]{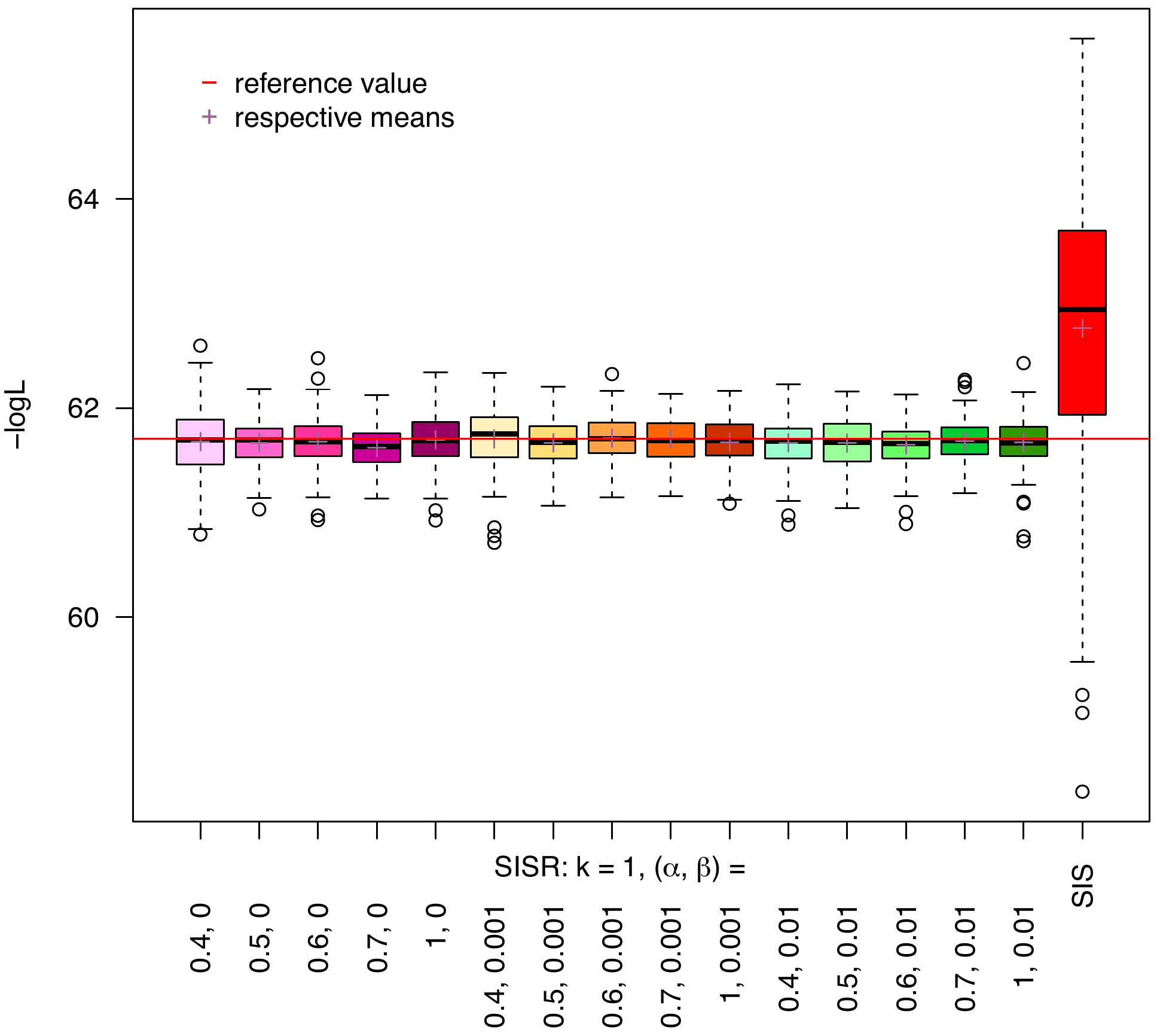}
\label{Boxplot_ncoalfixe_beta}}\\
\subfloat[]{
\includegraphics[scale=0.4]{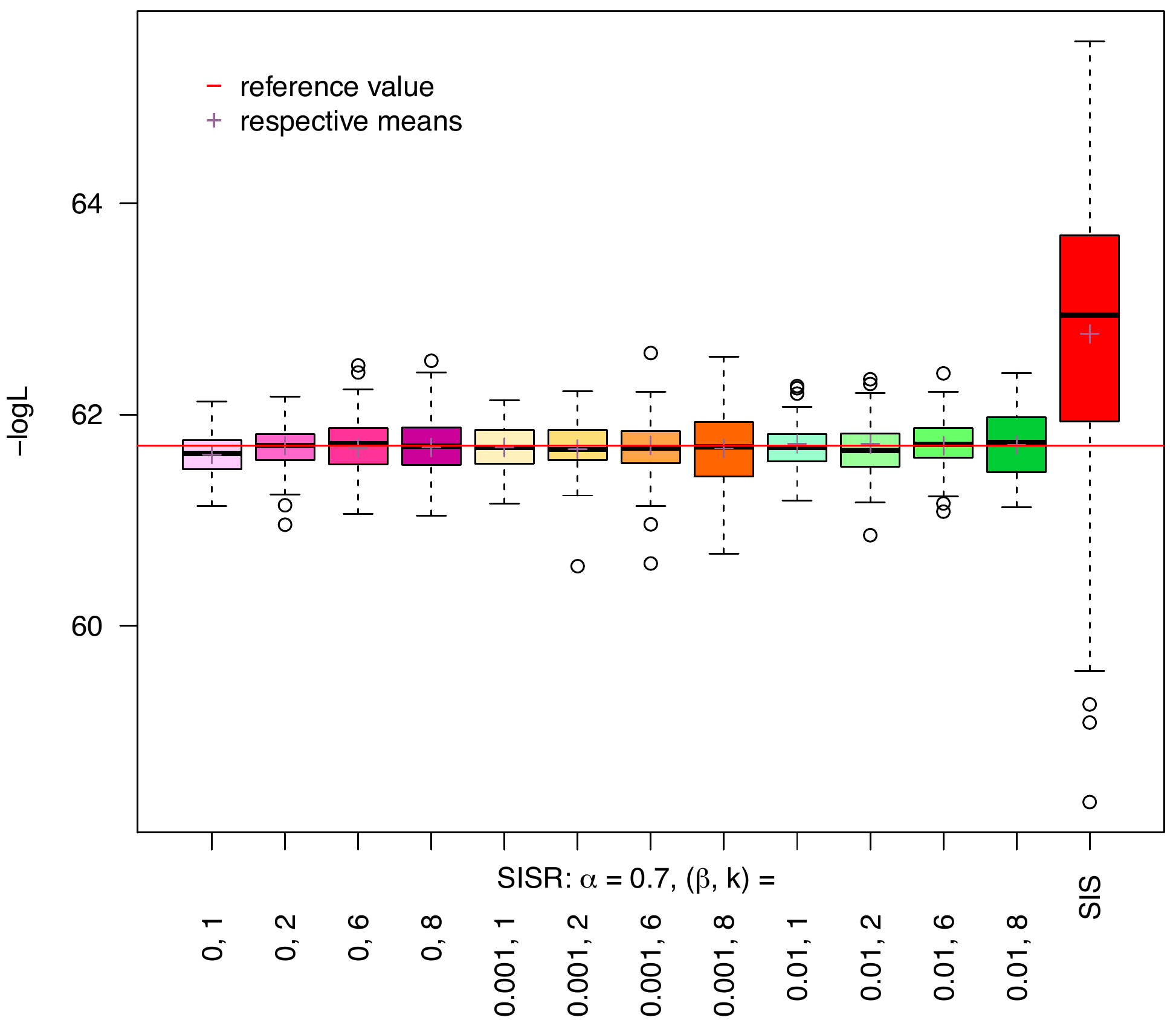}
\label{Boxplot_alphafixe_beta}}
\subfloat[]{
\includegraphics[scale=0.4]{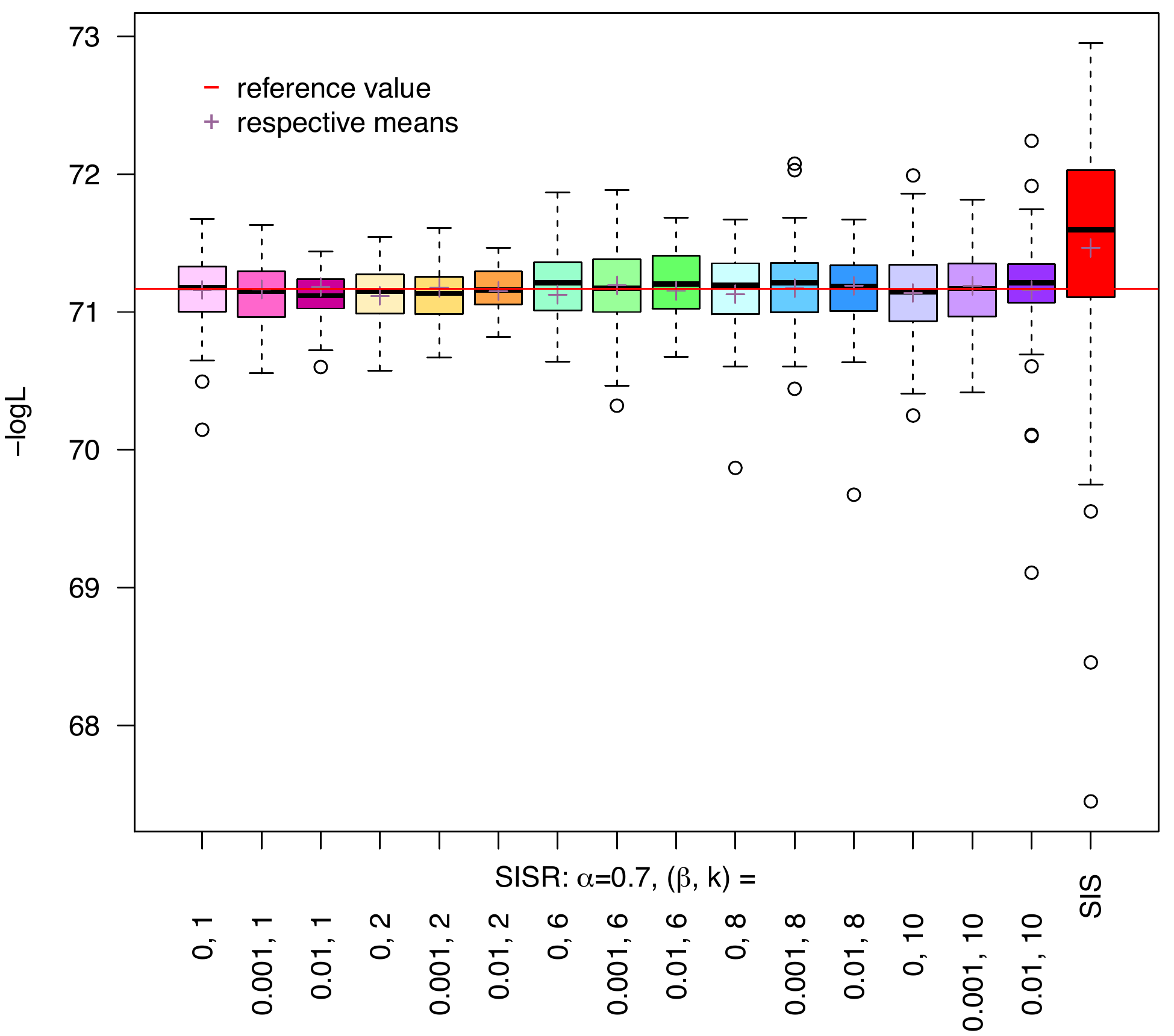}
\label{Boxplot_alphafixe_ncoal}}
\caption{\textbf{Boxplots of \(100\) estimates of the likelihood in a given parameter point with different inference algorithms.} The red box correspond to SIS estimates. The red horizontal line represents the reference value and the cross on each box represents the mean of the $100$ estimates of this box. Each block of three to five boxes with similar colors (pink, orange, green, bleu and mauve) corresponds to SISR with a fixed value: (a) of $\alpha$ for four different values of $k$ when $\beta=0.01$ is fixed for all the SISR estimations, (b) of $\beta$ for five different values of $\alpha$ when $k=1$ is fixed for all the SISR estimations, (c) of $\beta$ for four different values of $k$ when $\alpha=0.7$ is fixed for all the SISR estimations, and (d) of $k$ for three different values of $\beta$ when $\alpha=0.7$ is fixed for all the SISR estimations. See main text for details about $k$, $\alpha$ and $\beta$.
}
\label{Boxplots}
\end{figure}


\begin{figure}[!h]
\centering
\subfloat[]{
\includegraphics[scale=0.41]{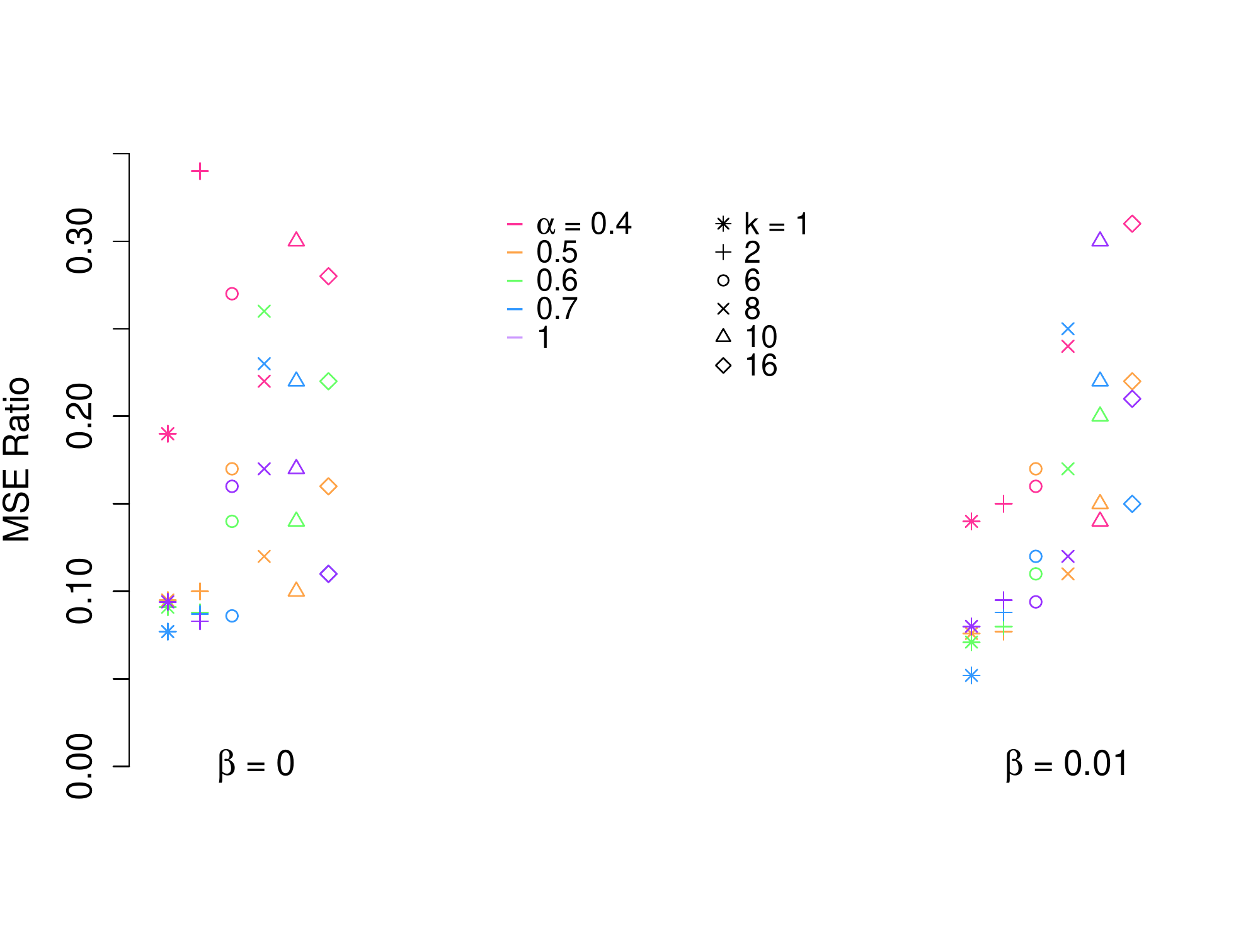}\label{0.4_0_0.25_400_RMSEfctbeta_decalagencoal}}
\subfloat[]{
\includegraphics[scale=0.41]{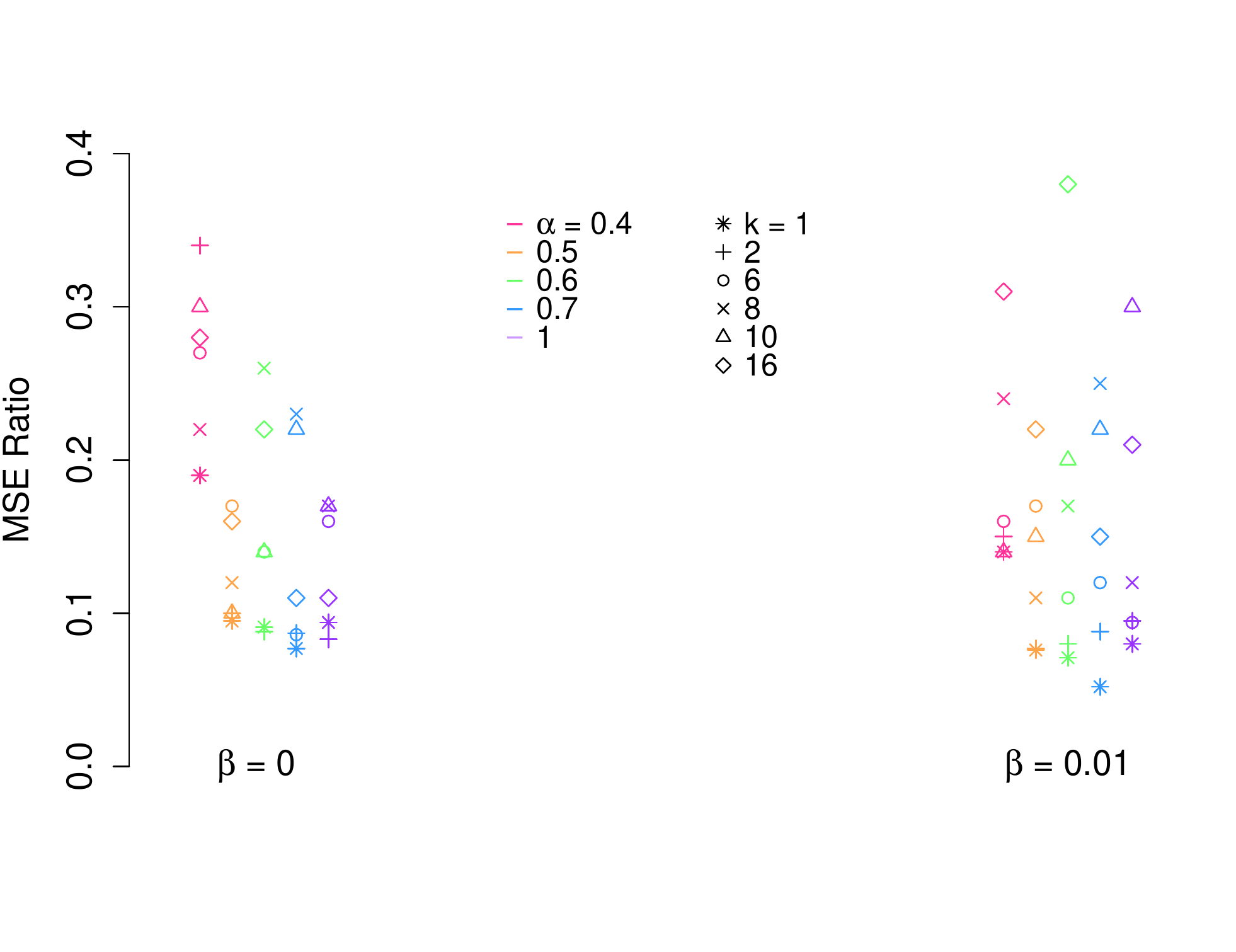}\label{0.4_0_0.25_400_RMSEfctbeta_decalagealpha}}\\
\subfloat[]{
\includegraphics[scale=0.7]{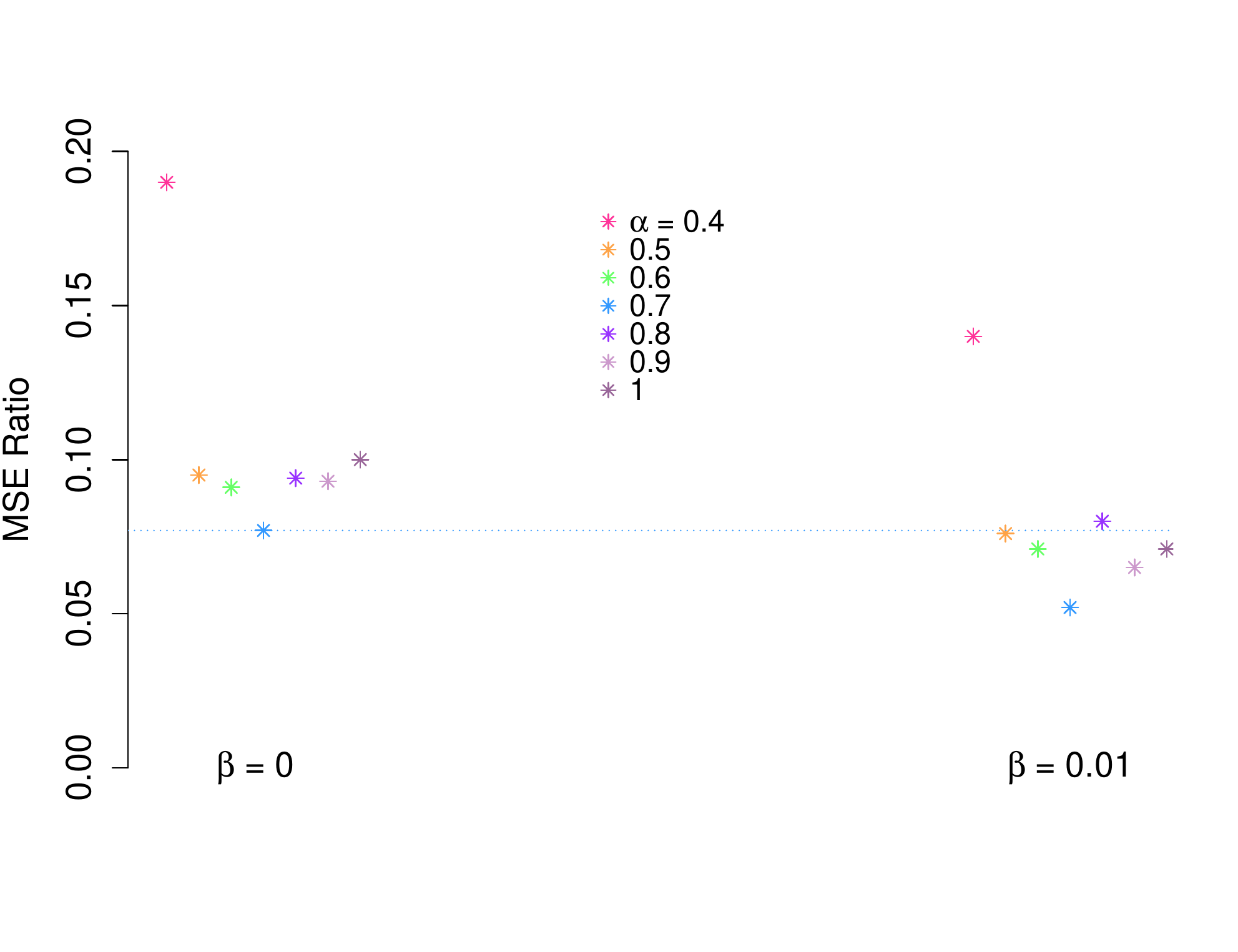}\label{0.4_0_0.25_400_RMSEfctbeta_1coal}}\\
\caption{
\textbf{MSE Ratios obtained with $\beta=0$ and $0.01$}.
(a) and (b) both represent the same MSE ratios but differ on the arrangement of the points. Each color corresponds to a value of $\alpha$ and each shape corresponds to a value of $k$, 
(a) each vertical alignment of five points corresponds to a fixed value of $k$ for five different values of $\alpha$, (b) each vertical alignment of six points corresponds to a fixed value of $k$ for six different values of $\alpha$.
(c) MSE Ratios obtained with $k=1$ and different values of $\alpha$. The horizontal blue line represents the lower MSE Ratio obtained with $k=1$ and $\beta=0$ among different values of $\alpha$.  See main text for details about $k$, $\alpha$ and $\beta$.}
\label{0.4_0_0.25_400_RMSE}
\end{figure}


\begin{figure}[!h]
\centering
\subfloat[\textbf{SIS} \(n_H=50\)]{%
\includegraphics[scale=0.65]{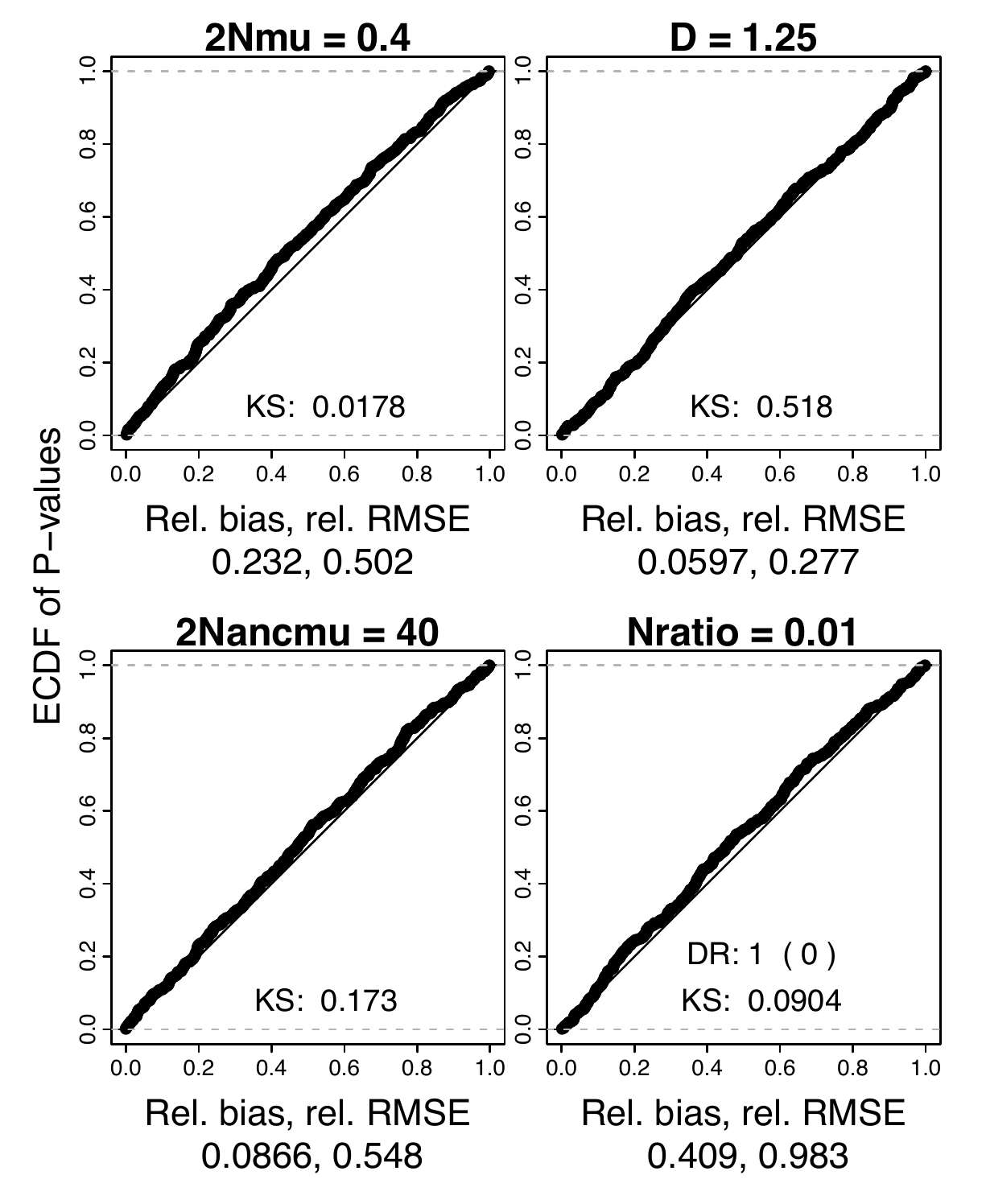}\label{0.4_0_1.25_40_50A_SIS}}
\subfloat[\textbf{SISR} \(n_H=50\)]{%
\includegraphics[scale=0.65]{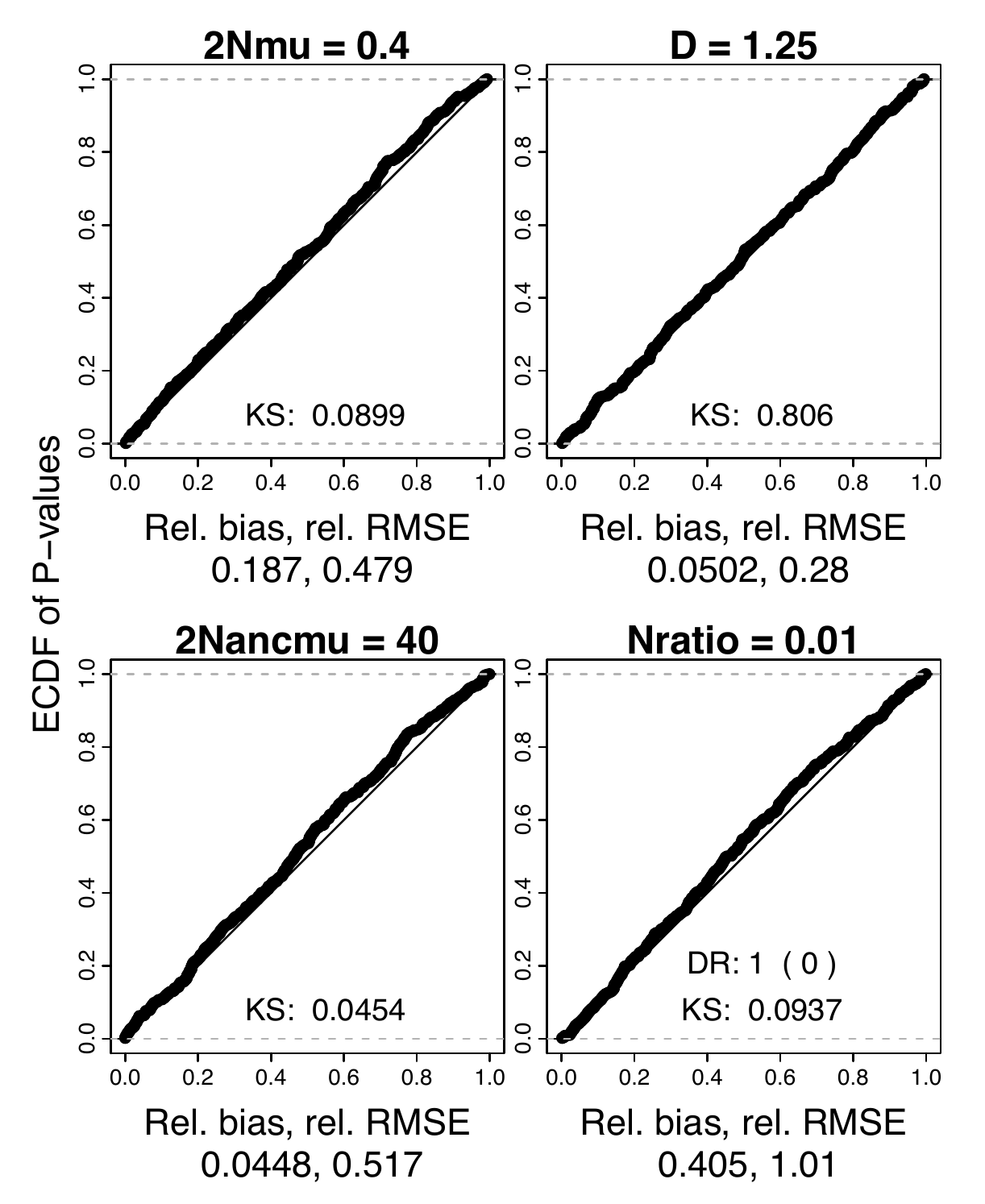}
\label{0.4_0_1.25_40_50A_SISR}} \\
\subfloat[\textbf{SIS} \(n_H=100\)]{%
 \includegraphics[scale=0.65]{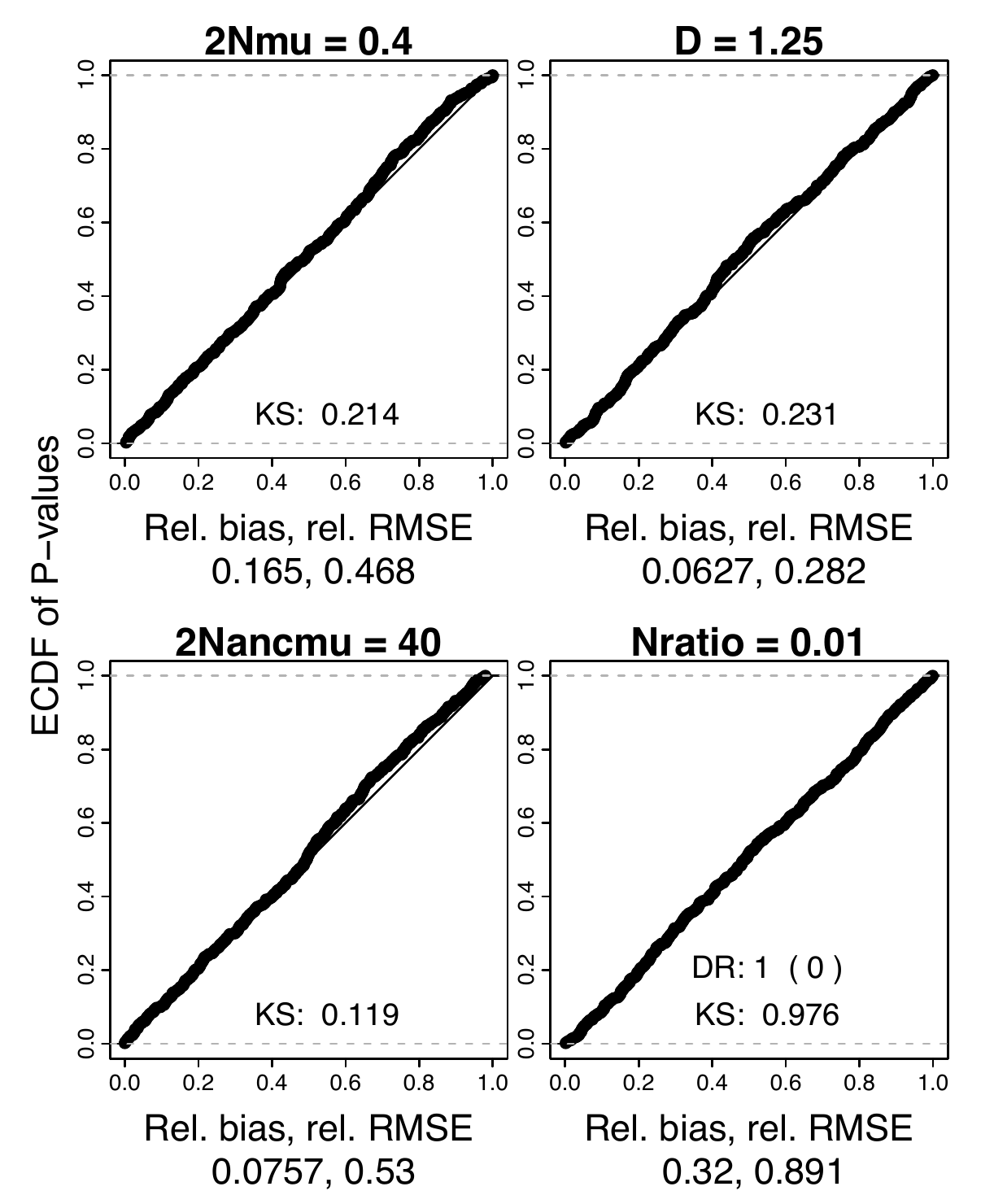}
\label{0.4_0_1.25_40_100A_SIS}}
\subfloat[\textbf{SISR} \(n_H=100\)]{%
\includegraphics[scale=0.65]{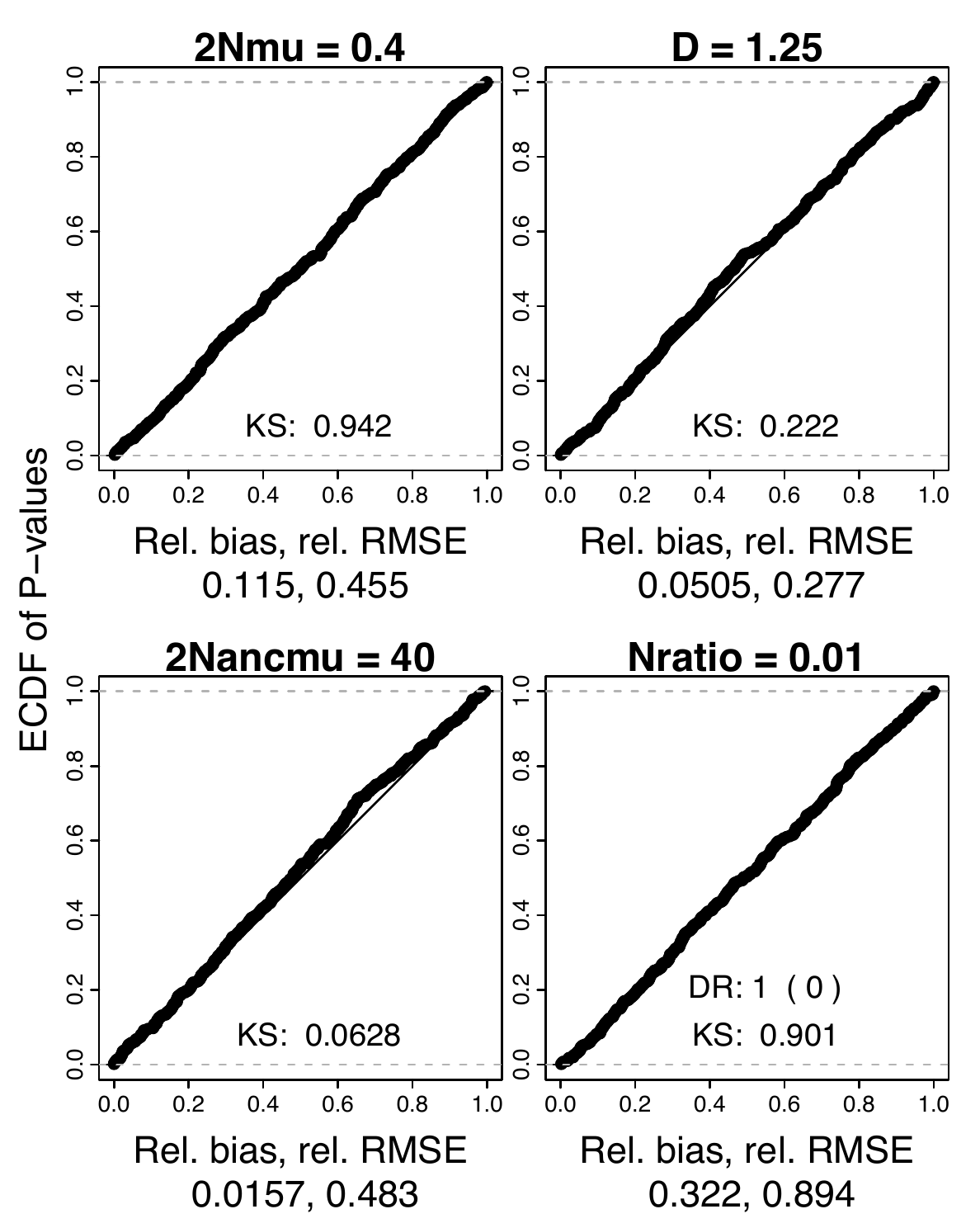}
\label{0.4_0_1.25_40_100A_SISR}}
\caption{\textbf{Empirical Cumulative Distribution Functions (ECDF) of p-values of Likelihood ratio tests for the scenario $\theta= 0.4$, $D =1.25$ and $\theta_\text{anc} = 40$.} Inference (a) with the SIS procedure with $n_H=50$ sampled histories (b) with the SISR procedure with $n_H=50$   (c)  SIS with $n_H=100$ (d) SISR with $n_H=100$, on  $500$ simulated data sets.
Relative bias and relative RMSE are also reported, and KS indicate the p-value of the Kolmogorov-Smirnov test for departure of LRT p-values distributions from uniformity.
}
\label{0.4_0_1.25_40_50A_100A}
\end{figure}


\begin{figure}[!h]
\centering
\subfloat[\textbf{SIS} \(n_H=100\)]{
\includegraphics[scale=0.47]{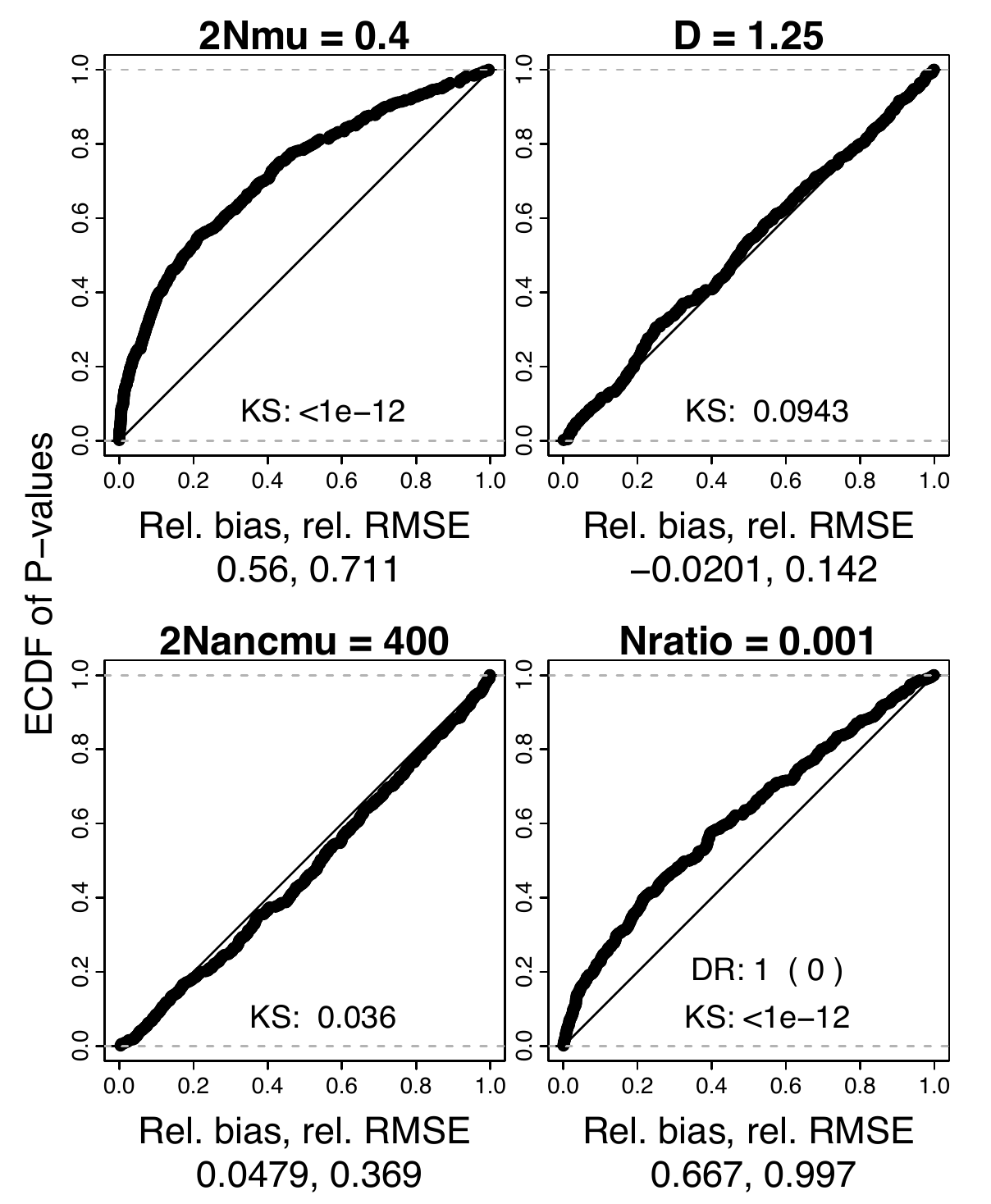}\label{0.4_0_1.25_400_100A_SIS}}
\subfloat[\textbf{SISR} \(n_H=100\)]{
\includegraphics[scale=0.47]{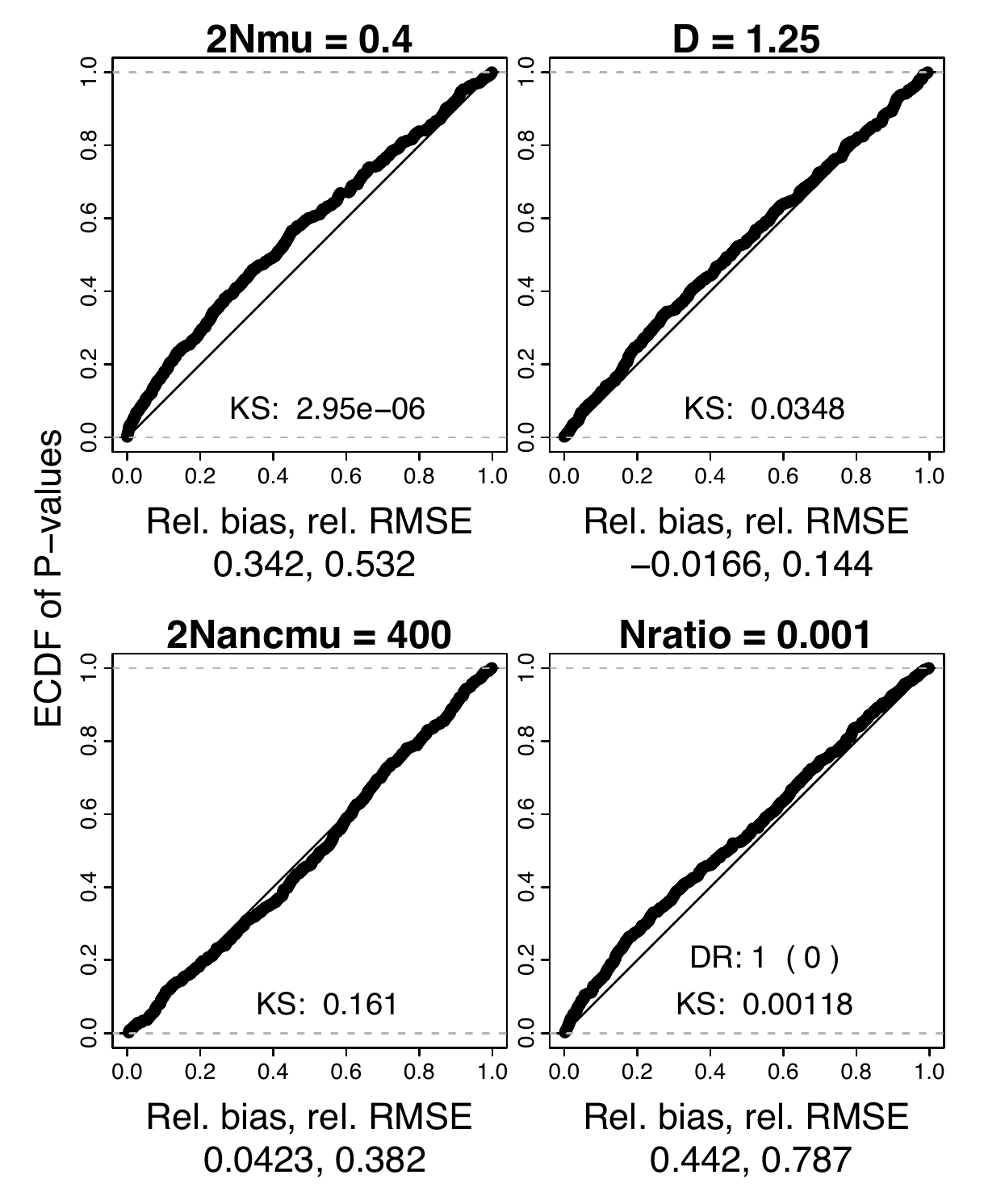}
\label{0.4_0_1.25_400_100A_SISR}} \\
\subfloat[\textbf{SIS} \(n_H=200\)]{
\includegraphics[scale=0.47]{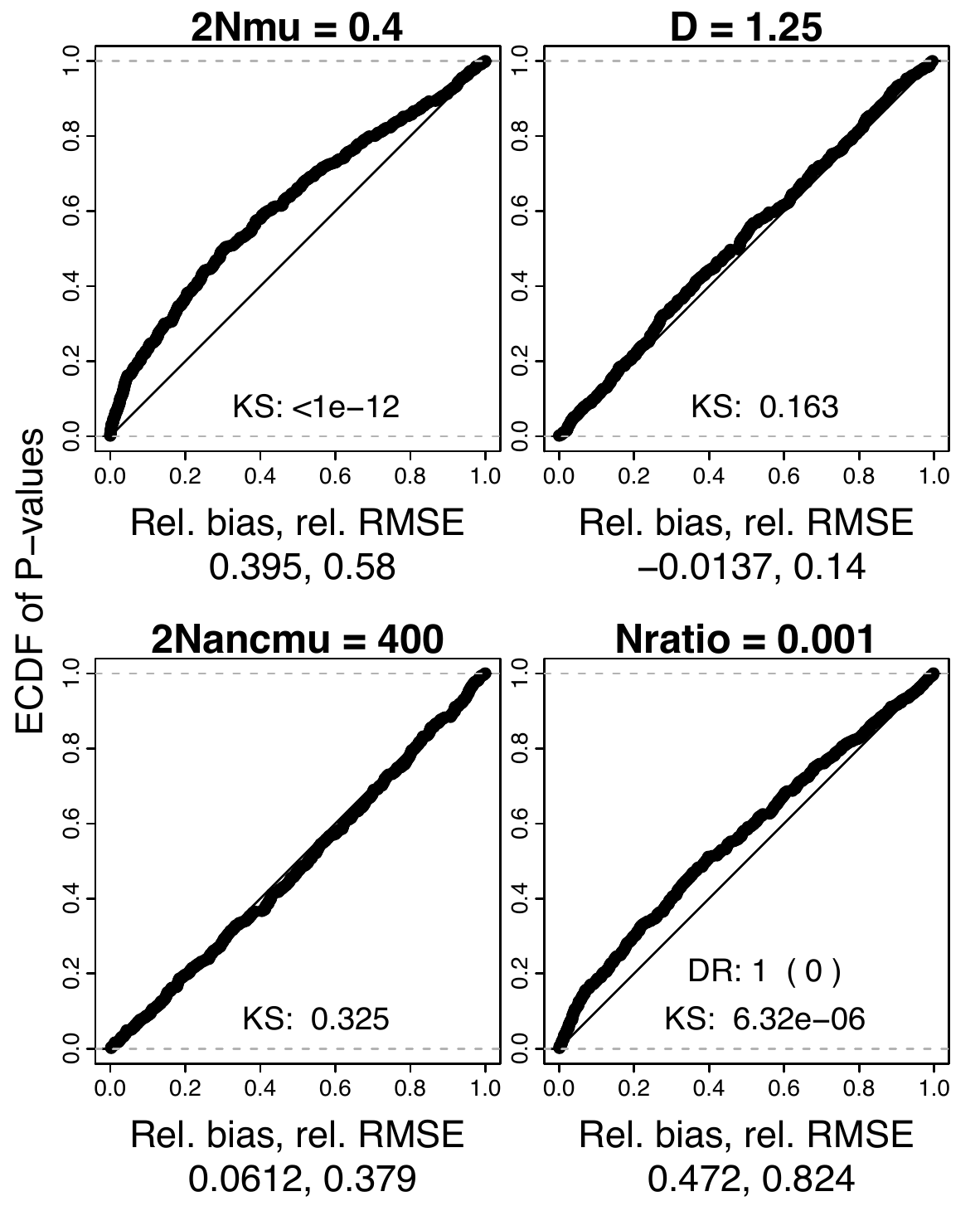}\label{0.4_0_1.25_400_200A_SIS}}
\subfloat[\textbf{SISR} \(n_H=200\)]{
\includegraphics[scale=0.47]{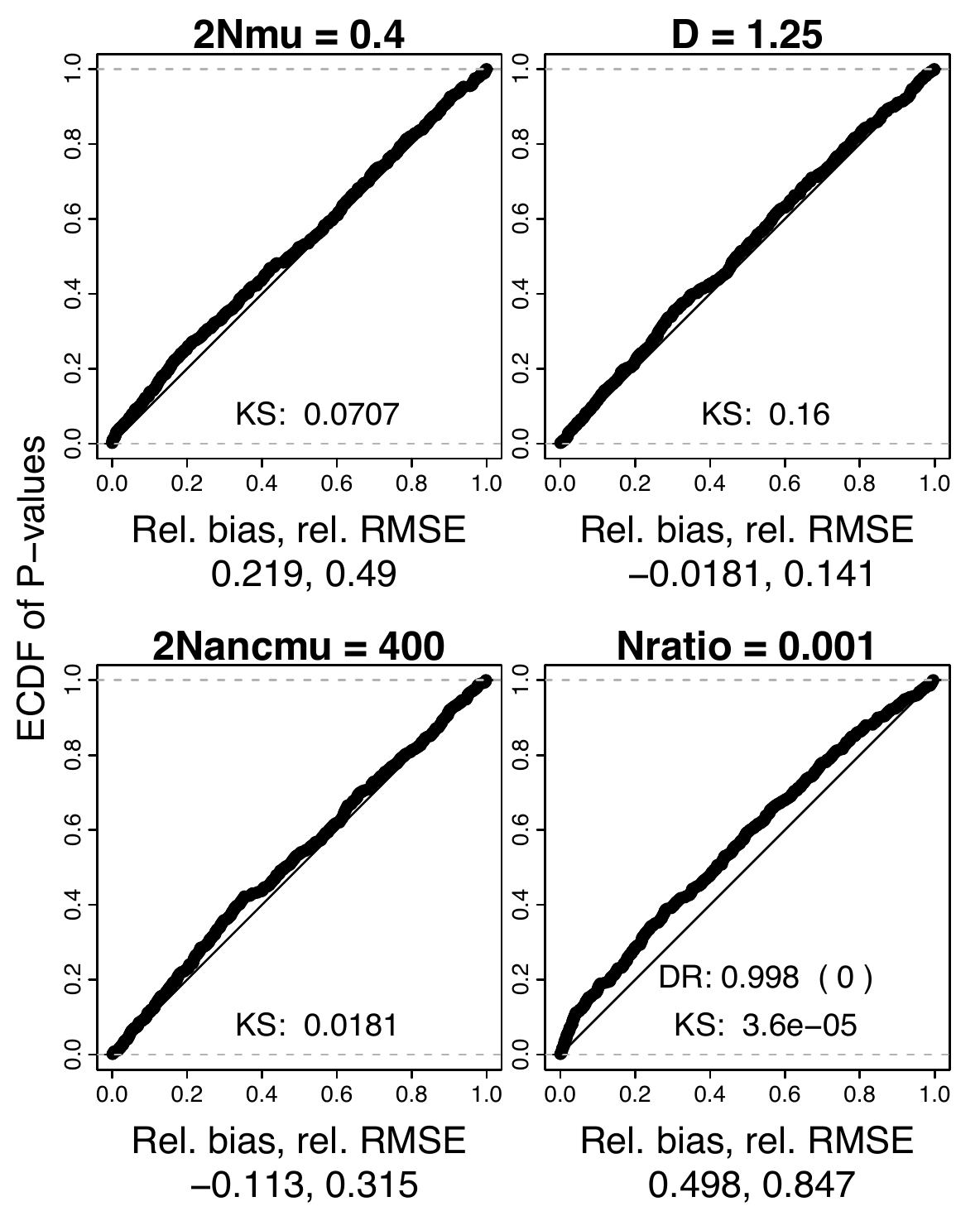}
\label{0.4_0_1.25_400_200A_SISR}} \\
\subfloat[\textbf{SIS}  \(n_H=400\)]{
 \includegraphics[scale=0.47]{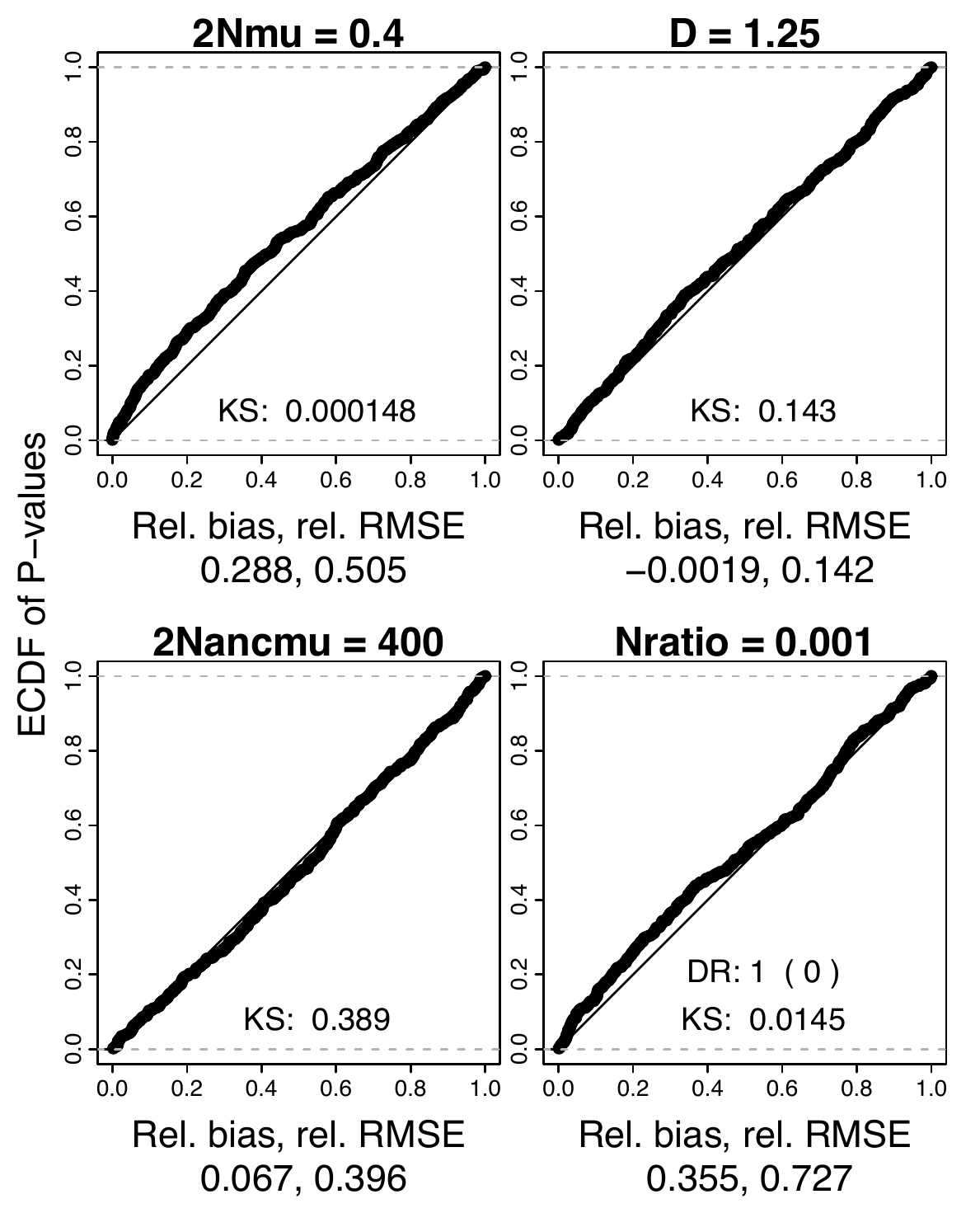}
\label{0.4_0_1.25_400_400A_SIS}}
\subfloat[\textbf{SISR}  \(n_H=400\)]{
\includegraphics[scale=0.47]{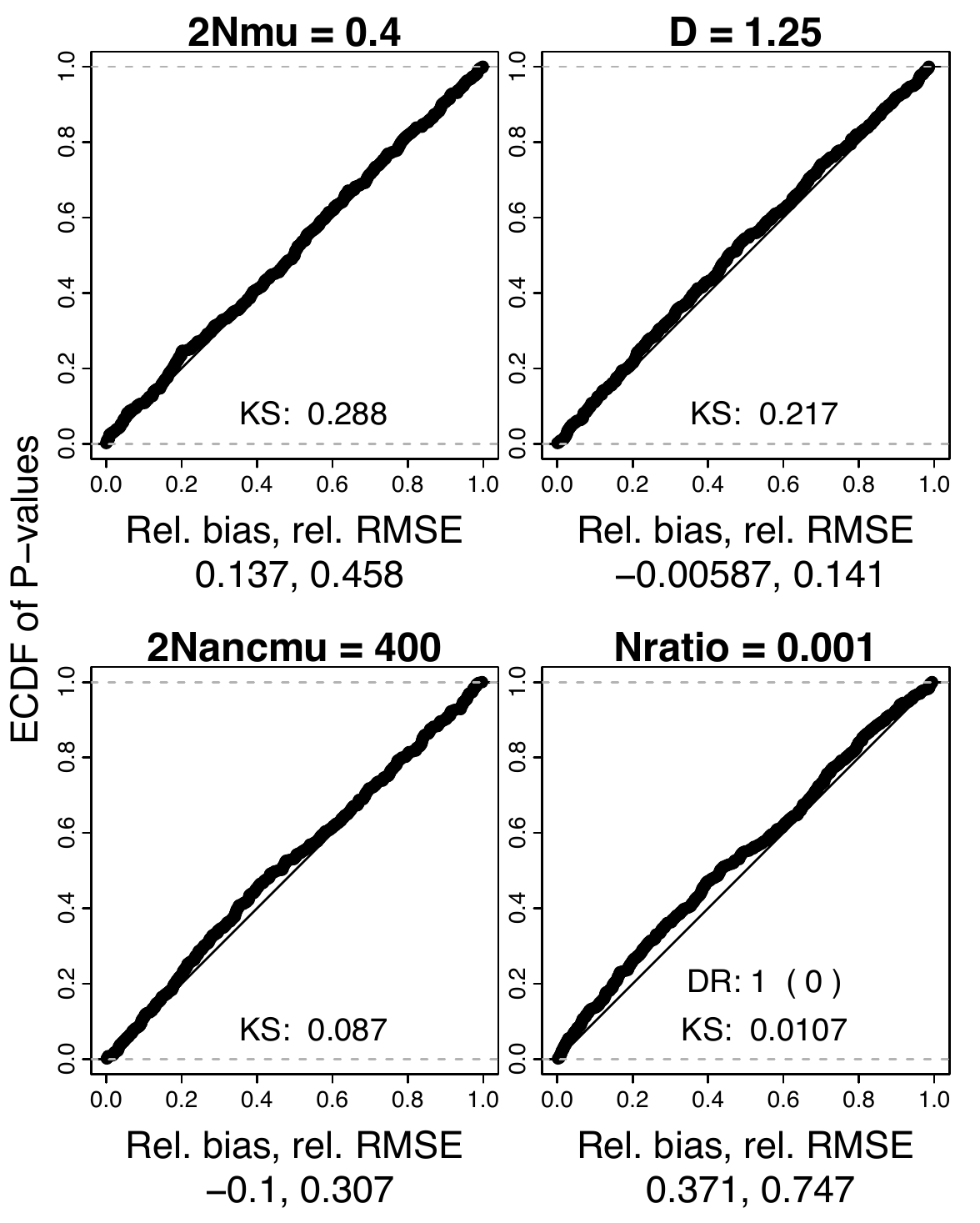}
\label{0.4_0_1.25_400_400A_SISR}}
\caption{\textbf{ECDF of p-values of Likelihood ratio tests for the scenario $\theta= 0.4$, $D =1.25$ and $\theta_\text{anc} = 400$.} (a) and (b) with $n_H=100$ (c) and (d) with $n_H=200$ and (e) and (f) with $n_H=400$, on $500$ data sets. See Fig.~\ref{0.4_0_1.25_40_50A_100A} for details.}
\label{0.4_0_1.25_400}
\end{figure}


\begin{figure}[!h]
\centering
\subfloat[\textbf{SIS} $n_H=2000$]{%
\includegraphics[scale=0.65]{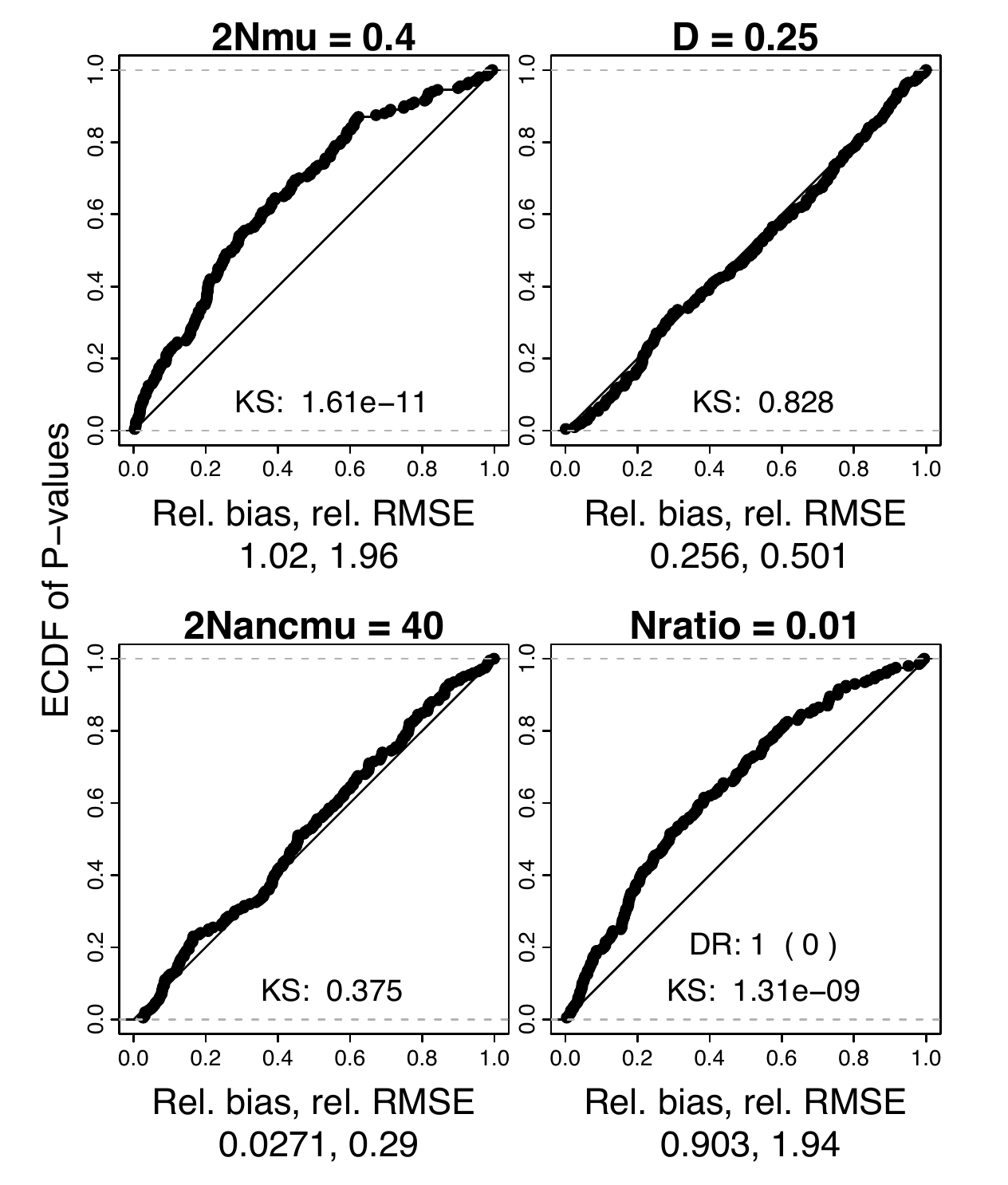}\label{0.4_0_0.25_40_2000A_SIS}}
\subfloat[\textbf{SISR} $n_H=2000$]{%
\includegraphics[scale=0.65]{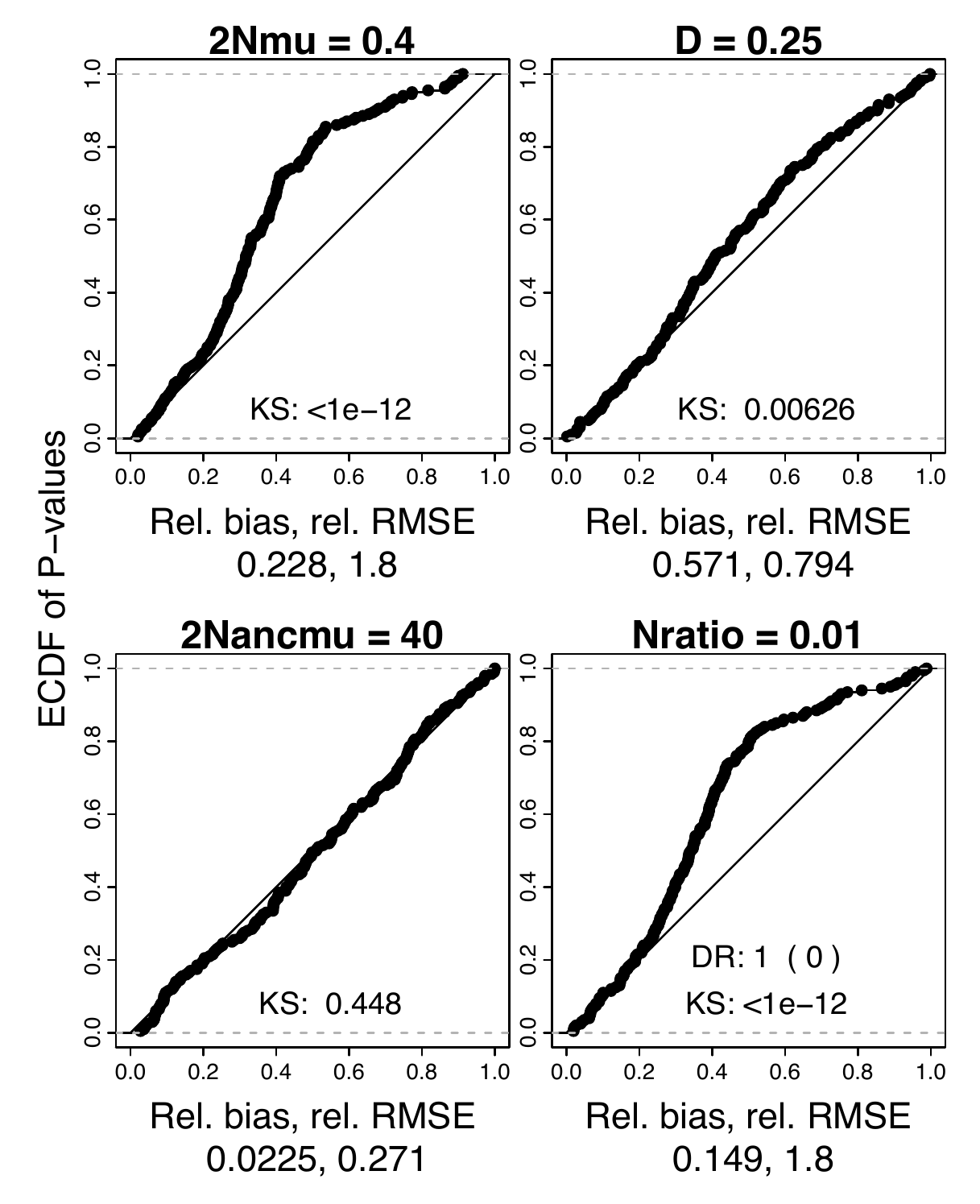}
\label{0.4_0_0.25_40_2000A_SISR}}
\caption{\textbf{ECDF of p-values of Likelihood ratio tests for the scenario $\theta= 0.4$, $D =0.25$ and $\theta_\text{anc} = 40$,}  with $n_H=2000$ sampled histories, on $200$ simulated data sets. See Fig.~\ref{0.4_0_1.25_40_50A_100A} for details.}
\label{0.4_0_0.25_40_2000A}
\end{figure}

\FloatBarrier
\newpage
 
\begin{figure}[!h]
\centering
\subfloat[\textbf{SIS}  \(n_H=2000\)]{
\includegraphics[scale=0.65]{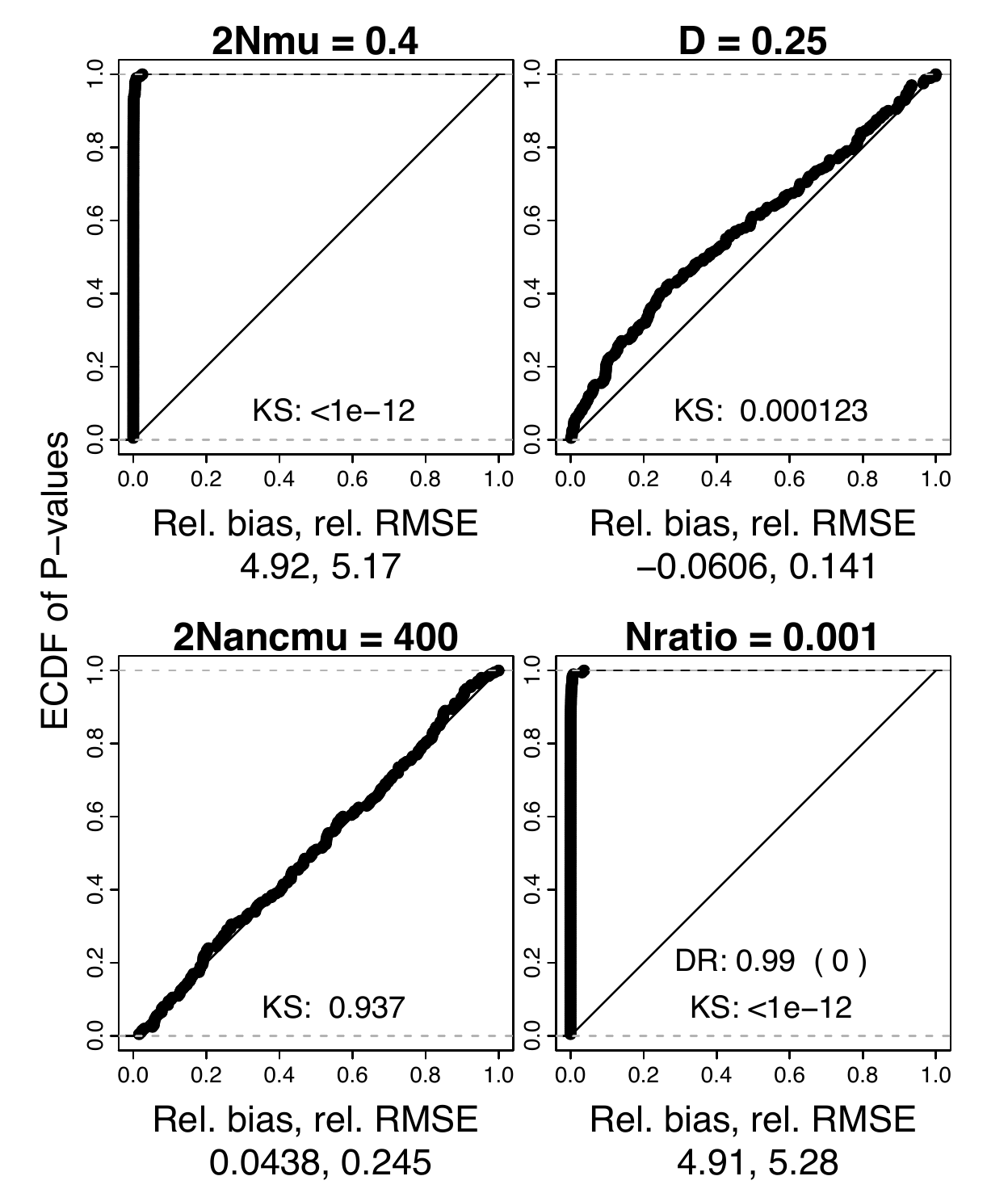}\label{0.4_0_0.25_400_2000A_SIS}}
\subfloat[\textbf{SIS}  \(n_H=20000\)]{
\includegraphics[scale=0.65]{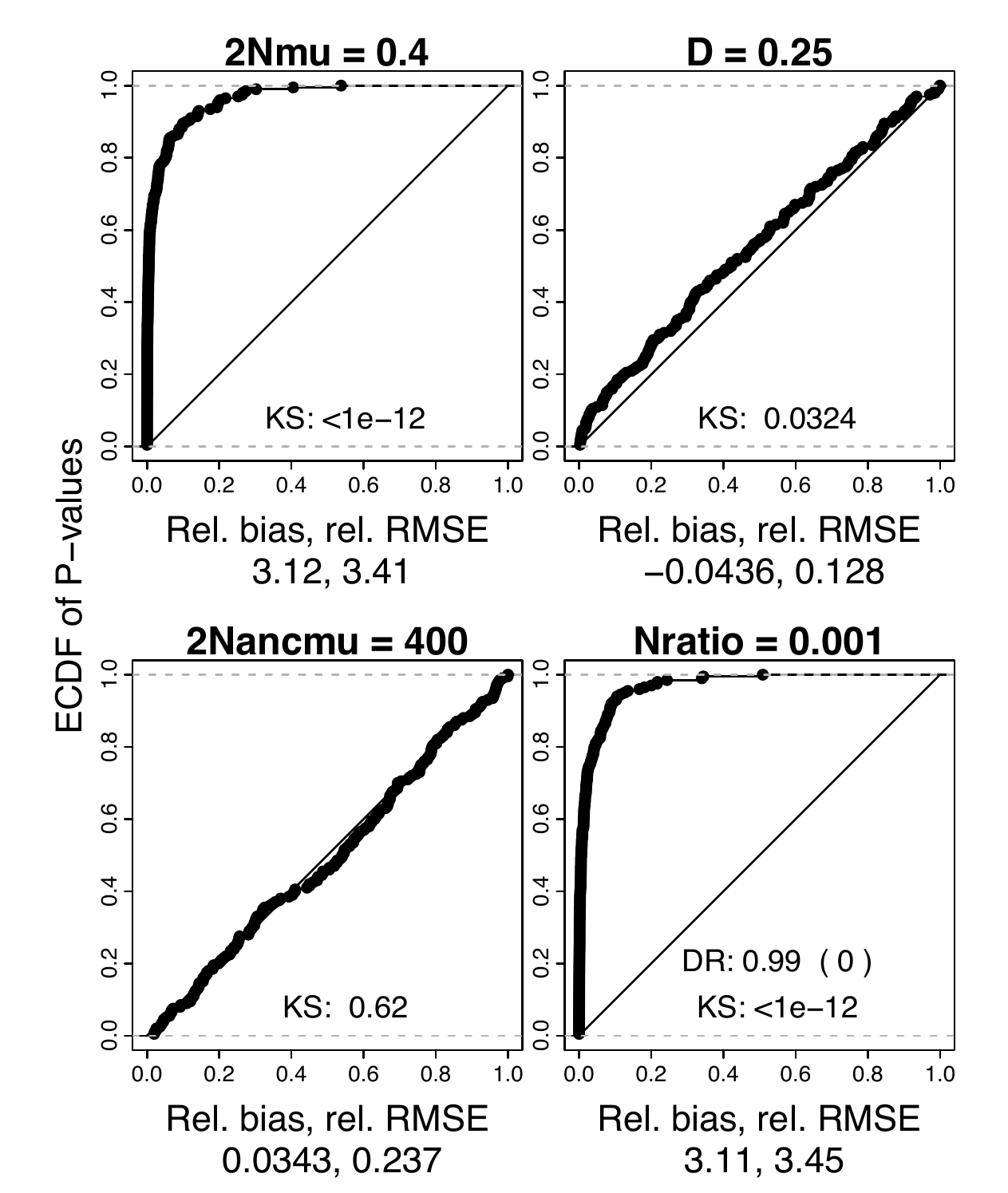}
\label{0.4_0_0.25_400_20000A_SIS}} \\
\subfloat[\textbf{SIS}  \(n_H=200000\)]{
\includegraphics[scale=0.65]{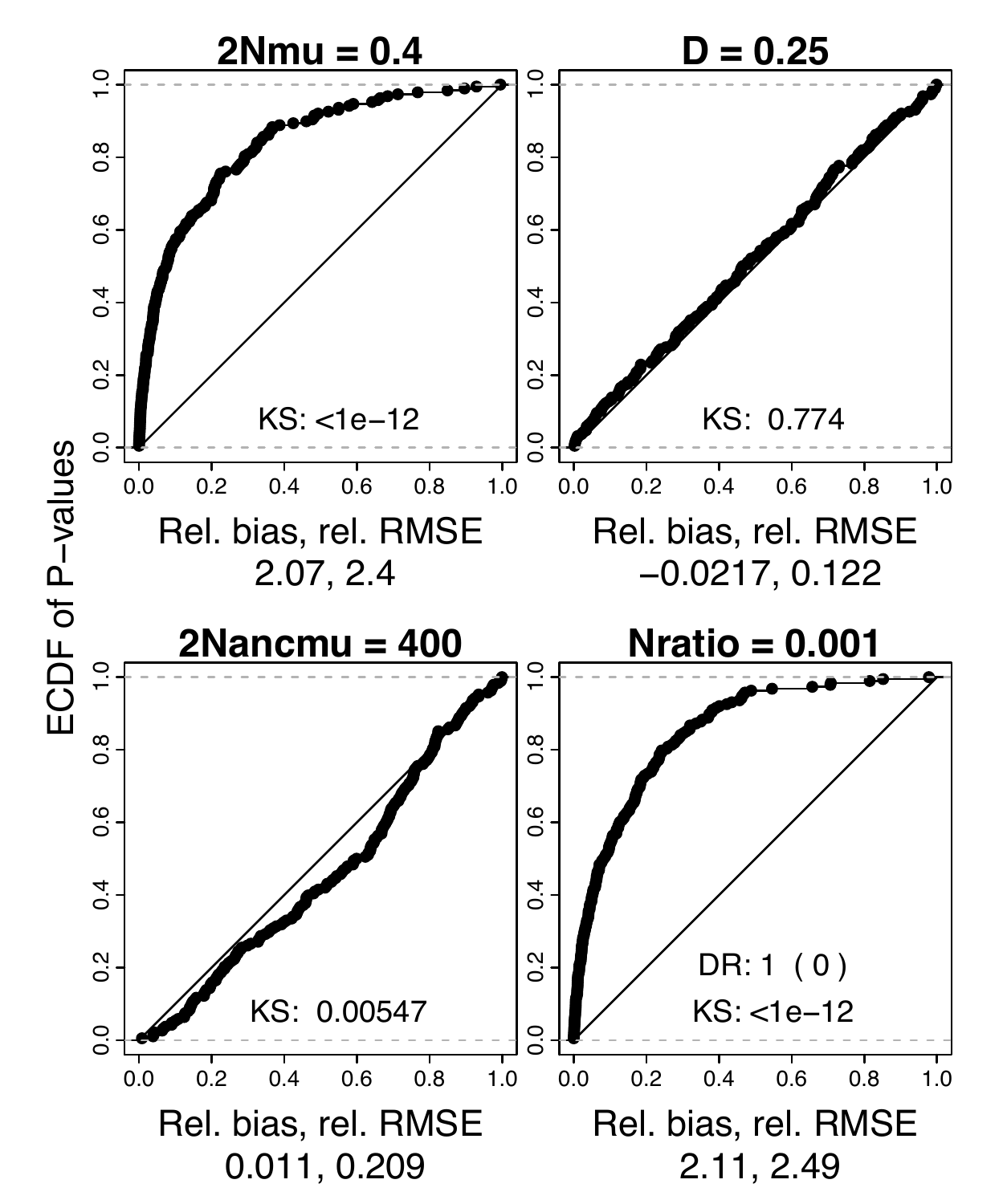}\label{0.4_0_0.25_400_200000A_SIS}}
\subfloat[\textbf{SISR} \(n_H=2000\)]{
\includegraphics[scale=0.65]{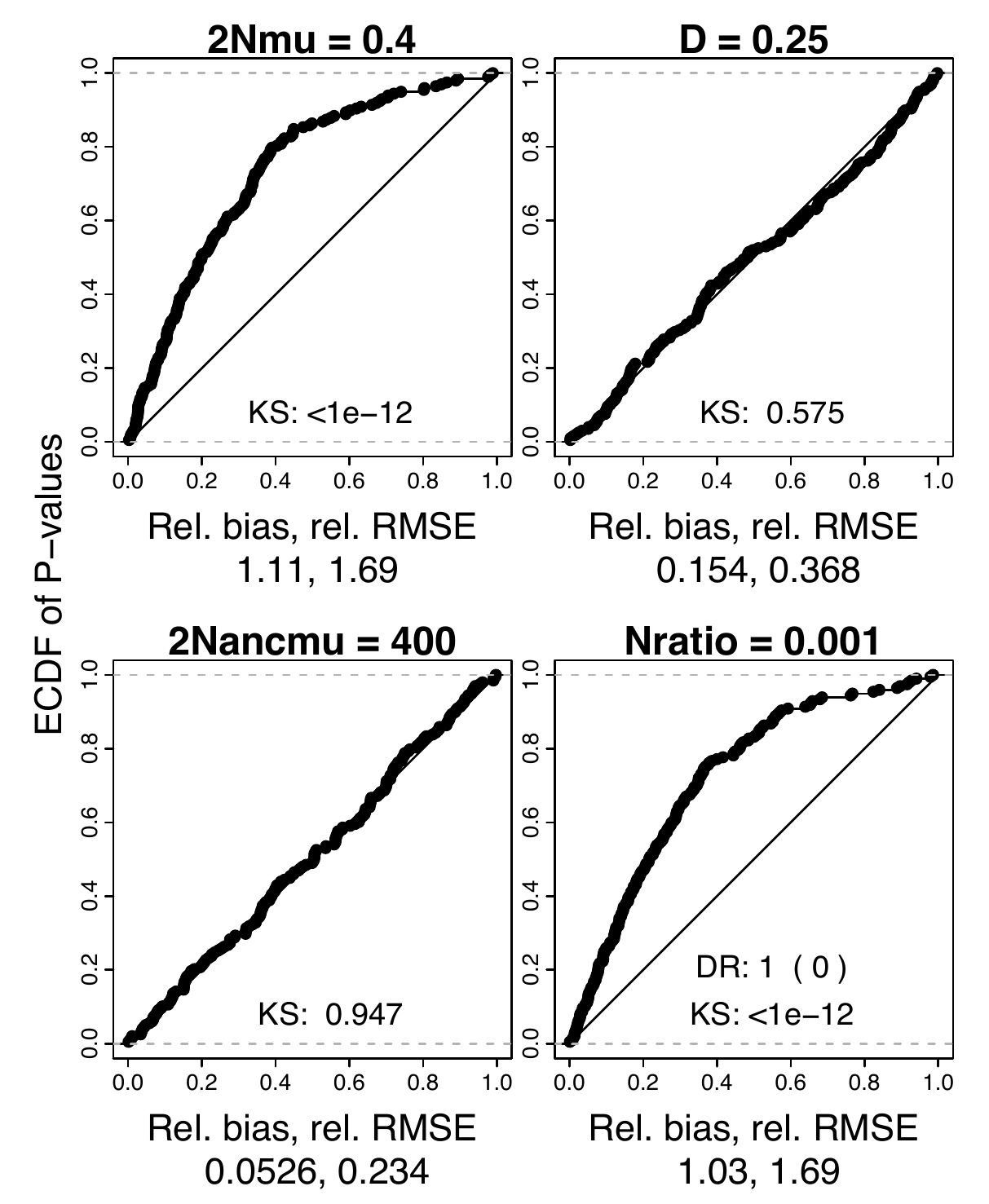}
\label{0.4_0_0.25_400_2000A_SISR}}
\caption{\textbf{ECDF of p-values of Likelihood ratio tests for the scenario $\theta= 0.4$, $D =0.25$ and $\theta_\text{anc} = 400$.} (a) SIS with $n_H=2000$ sampled histories (b) SIS with $n_H=20000$ sampled histories (c) SIS with $n_H=200000$ sampled histories (d) SISR with $n_H=2000$ sampled histories, on $200$ simulated data sets. See Fig.~\ref{0.4_0_1.25_40_50A_100A} for details.}
\label{0.4_0_0.25_400}
\end{figure}

\newpage

\begin{figure}[!h]
\centering
\subfloat[\textbf{SIS}]{
\includegraphics[scale=0.9]{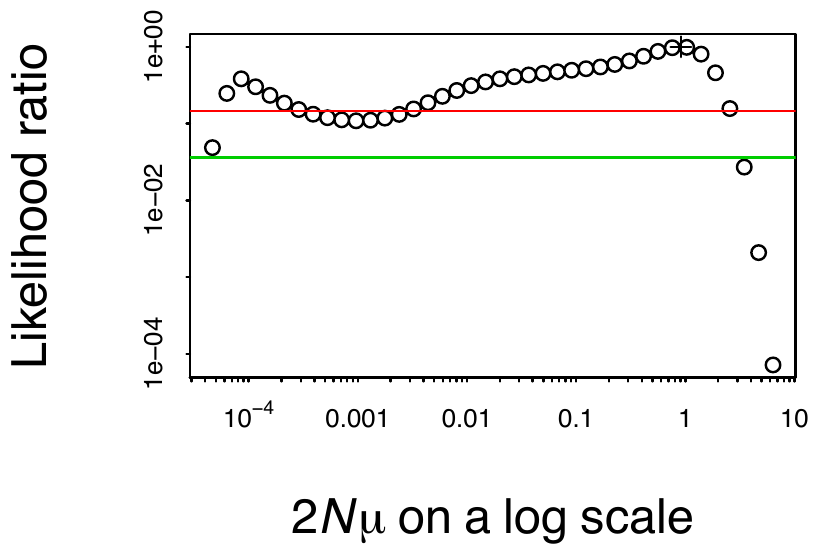}\label{0.4_0_0.25_40_2000A_profiles_SIS_Theta}}
\subfloat[\textbf{SISR}]{
\includegraphics[scale=0.9]{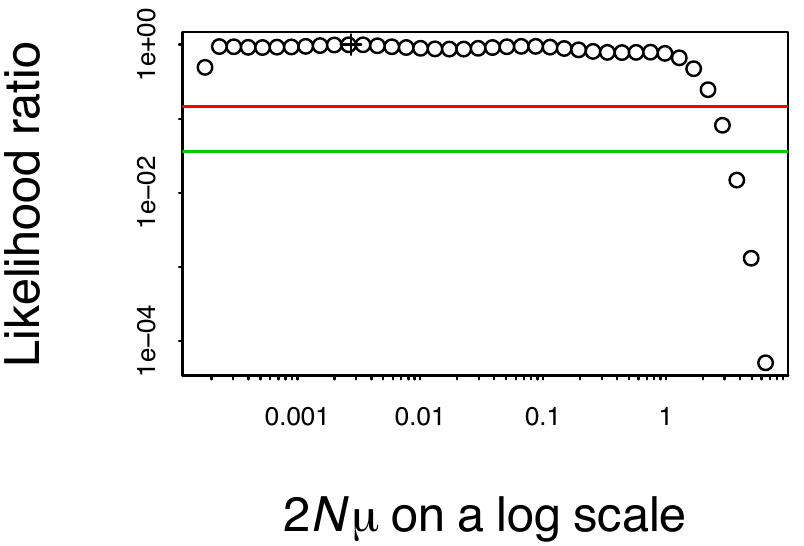}
\label{0.4_0_0.25_40_2000A_profiles_SISR_Theta}} \\
\subfloat[\textbf{SIS}]{
\includegraphics[scale=0.9]{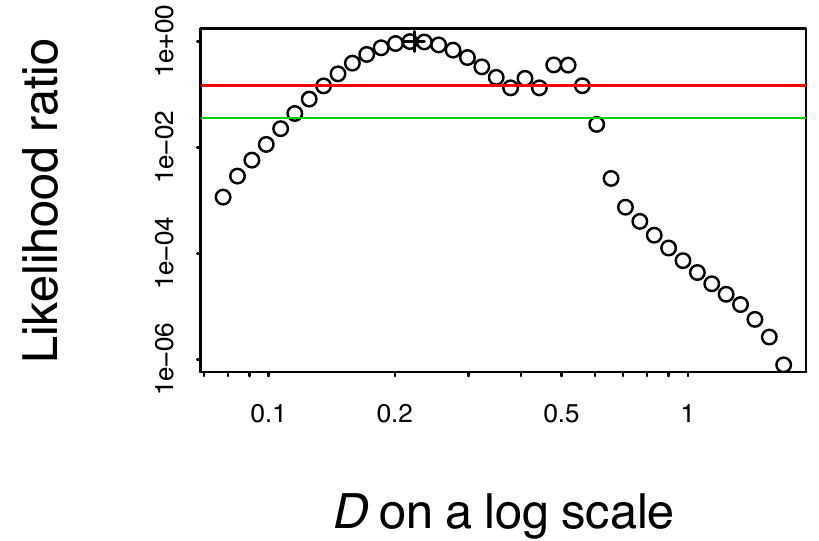}\label{0.4_0_0.25_40_2000A_profiles_SIS_D}}
\subfloat[\textbf{SISR}]{
\includegraphics[scale=0.9]{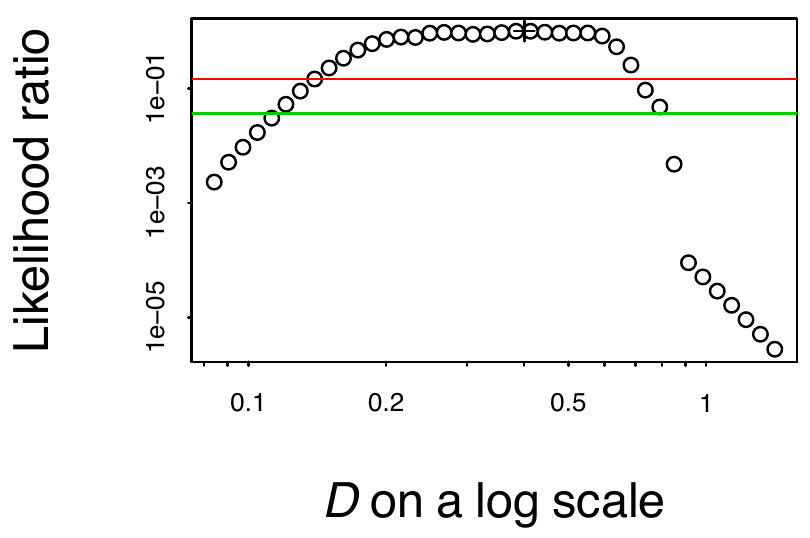}
\label{0.4_0_0.25_40_2000A_profiles_SISR_D}}\\
\subfloat[\textbf{SIS}]{
\includegraphics[scale=0.9]{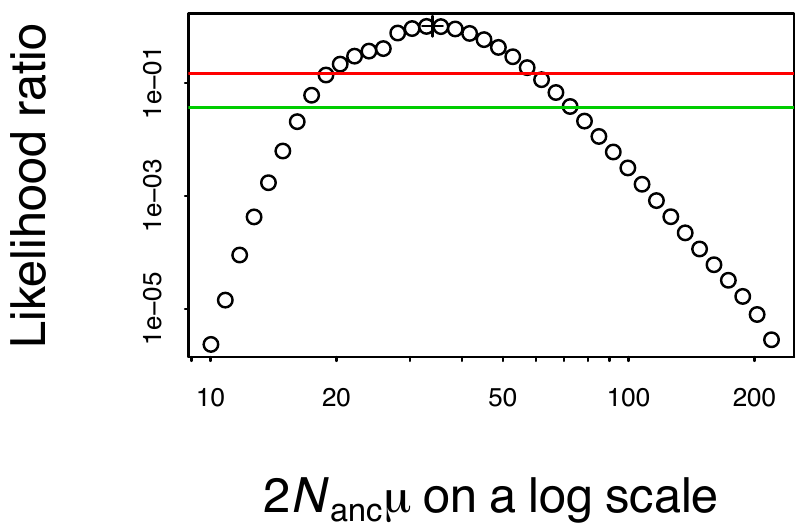}\label{0.4_0_0.25_40_2000A_profiles_SIS_ThetaAnc}}
\subfloat[\textbf{SISR}]{
\includegraphics[scale=0.9]{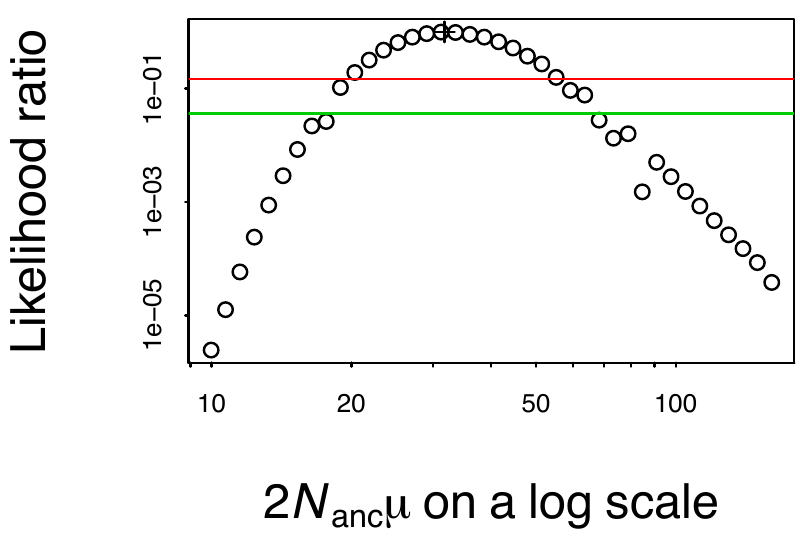}
\label{0.4_0_0.25_40_2000A_profiles_SISR_ThetaAnc}} \\
\caption{\textbf{One-parameter profile likelihood ratios}.
(a) and (b) represent the likelihood profile function divided by its maximum value $\theta \mapsto \hat{L}_{\theta}(\hat{D}_\theta,\hat{\theta}_{\text{anc }\theta})/\hat{L}(\hat{\theta},\hat{D}_\theta,\hat{\theta}_{\text{anc } \theta})$, (c) and (d) represent in the same way the profile likelihood ratio for $D$ and (e) and (f) represent the profile likelihood ratio for $\theta_\text{anc}$, estimated respectively with SIS (left) or with SISR (right). See Appendix.~\ref{ProfileLikelihood} for details on profile likelihood.}
\label{0.4_0_0.25_40_2000A_profiles}
\end{figure}

\end{document}